\pdfoutput=1
\documentclass{article}

\usepackage[utf8]{inputenc}
\usepackage{enumitem}
\usepackage{setspace}
\usepackage{mathtools}
\usepackage{tikz-cd}
\usepackage{lipsum}
\usepackage{tikz}
\usepackage{mathrsfs}
\usetikzlibrary{arrows}
\usepackage{amssymb}
\usetikzlibrary{decorations.pathmorphing}
\usepackage{bbold}
\usepackage{amsmath,amsthm}
\usepackage{xspace}
\usepackage{faktor}
\usepackage{calrsfs}
\usepackage{verbatim}
\usepackage{quiver}
\usepackage[pdfa]{hyperref}

\newtheorem{theorem}{Theorem}[section]
\newtheorem{proposition}[theorem]{Proposition}
\newtheorem{lemma}[theorem]{Lemma}

\newtheorem{corollary}[theorem]{Corollary}

\newtheorem{definition}[theorem]{Definition}

\newtheorem{remark}[theorem]{Remark}

\newcommand{\mv}{\mathbb{MV}}
\newcommand{\cat}[1]{\mathbb{#1}}
\newcommand{\mono}{\rightarrowtail}
\newcommand{\inv}[1]{#1^{-1}}
\newcommand{\pol}[1]{#1^{\perp}}

\newcommand{\pt}[2]{Pt_{#1}{#2}}

\newcommand{\lab}{\ell \mathbb{A}b}
\newcommand{\ulab}{u\ell \mathbb{A}b}

\newcommand{\Extmv}{\Ext \mv}
\newcommand{\CExtmv}{\CExt_{s\mv}\mv}
\newcommand{\Extdmv}{\Ext^2 \mv}
\newcommand{\ExtCExtmv}{\Ext\CExt_{s\mv}\mv}
\newcommand\myfunc[5]{\begin{aligned}
  #1 \colon #2 &\to #3\\
  #4 &\mapsto #5
\end{aligned}}

\newcommand{\ibang}[1]{\iota_{#1}}
\newcommand{\pre}[1]{\mathcal{#1}}

\newcommand{\facte}{\mathcal{E}}
\newcommand{\factm}{\mathcal{M}}

\newcommand{\zeros}{\mathcal{Z}}
\newcommand{\zeroarrows}{N_{\zeros}}

\DeclareMathOperator{\kerr}{ker}

\DeclareMathOperator{\Eq}{Eq}

\DeclareMathOperator{\Ideals}{Ideals}

\DeclareMathOperator{\Rad}{Rad}
\DeclareMathOperator{\Inf}{Inf}

\DeclareMathOperator{\Ext}{Ext}
\DeclareMathOperator{\CExt}{CExt}

\DeclareMathOperator{\Arr}{Arr}

\DeclareMathOperator{\torsione}{(\mathcal{T},\mathcal{F})}

\DeclareMathOperator{\terminal}{\boldsymbol{1}}
\DeclareMathOperator{\initial}{\boldsymbol{2}}
\usepackage[english]{babel}	

\title{A Galois Theory and a Pretorsion Theory in MV-Algebras}
\author{Andrea Cappelletti}

\begin{document}

\maketitle

\begin{abstract}
    We explore the relationship between the category of MV-algebras and its full subcategories of perfect and semisimple algebras, showing that this pair of subcategories defines a pretorsion theory. We study the Galois structure associated with the reflection of semisimple MV-algebras, proving that it is admissible from the point of view of categorical Galois theory and characterizing the corresponding central extensions. These central extensions are themselves reflective into the category of surjective MV-algebra morphisms, and the corresponding adjunction is admissible, too. Thanks to this observation, we characterize higher central extensions, and we use them to define a non-pointed version of commutators between ideal subalgebras.
\end{abstract}

\section{Introduction}

MV-algebras \cite{chang1958algebraic} are a class of algebraic structures used in mathematics to study non-classical logic. They generalize the notion of Boolean algebras, which are used to model classical logic, and can be seen as a mathematical system for reasoning under uncertainty.  MV-algebras have applications in various areas of mathematics, such as logic, algebra, topology, and computer science, and can be used to model fuzzy logic, intuitionistic logic, and quantum logic.\\
More specifically, an MV-algebra $A$ is a set equipped with an operation $\oplus$, which is both associative, commutative, and has a neutral element $0$, and an operation $\lnot$, such that the following properties hold: $\lnot \lnot x = x$, $x \oplus \lnot 0 = \lnot 0$, and
$\lnot (\lnot x \oplus y) \oplus y = \lnot (\lnot y \oplus x) \oplus x$.
These conditions are intended to capture some properties of the real unit interval $[0, 1]$ equipped with negation $\lnot x = 1 - x$ and truncated addition $x \oplus y = \min(1, x + y)$.\\
MV-algebras are a powerful tool to model Łukasiewicz calculus mathematically \cite{chang1958algebraic}. In fact, in Łukasiewicz calculus, which is a many-valued logic, the truth values are not limited to just ``true'' or ``false''; instead, they can take on any value from a continuous interval, typically $[0,1]$. In this sense, MV-algebras provide a natural framework for modeling Łukasiewicz calculus, as they can represent the truth values as elements in the interval $[0,1]$, and develop an algebraic structure that can naturally manipulate these values.\\
In this paper, we will investigate some categorical-algebraic properties of the variety of MV-algebras, denoted by $\mv$. Although $\mv$ is not a semi-abelian category, it possesses many important properties. Specifically, a category is semi-abelian \cite{semiabel} if it is a pointed finitely cocomplete category which is Barr-exact \cite{barr} and protomodular \cite{bourn1991normalization}. However, $\mv$ fails to satisfy the condition of having isomorphic initial and terminal objects; the initial object is the algebra $\{0,1\}$, while the terminal object is the algebra $\{0=1\}$. Nonetheless, being a variety of universal algebras, $\mv$ is complete, cocomplete, and Barr-exact. Furthermore, we will show that $\mv$ is also a protomodular category. To summarize, although $\mv$ is not a semi-abelian category, it still possesses many important algebraic-categorical properties, making it an interesting object of study in its own right.\\
In the second part of the paper, our study revolves around a Galois theory over $\mv$, induced by the subcategories of perfect algebras and semisimple algebras.
Given an MV-algebra $A$, we define its radical, denoted with $\Rad(A)$, as the intersection of all maximal ideals of $A$. It has been shown \cite{salas1980estudio} that $\Rad(A)$ consists precisely of those elements $a\in A$ that satisfy the inequality $na\leq \lnot a$ for every natural number $n$.
This notion of radical has several implications in the study of MV-algebras. In particular, it naturally leads to the definition of two classes of MV-algebras: the perfect MV-algebras and the semisimple MV-algebras.
An MV-algebra $A$ is said to be \emph{perfect} if it can be expressed as the union of its radical $\Rad(A)$ and the set $\lnot\Rad(A)$, which consists of all elements obtained as the negation of the elements in $\Rad(A)$. Interestingly, it has been shown (for a proof of this fact see \cite{cignoli2013algebraic}) that the category of non-trivial perfect MV-algebras is equivalent to that of lattice-ordered abelian groups.
An MV-algebra is said to be \emph{semisimple} if its radical is trivial. This notion of semisimplicity plays a crucial role in the study of the structure of MV-algebras, and it is intimately related to the notion of simplicity in other areas of algebra.\\
We denote by $p\mv$ the full subcategory of $\mv$ whose objects are the perfect MV-algebras, and by $s\mv$ the full subcategory of $\mv$ whose objects are the semisimple MV-algebras. Notably, we observe that $p\mv \cap s\mv = \{\terminal, \initial\}$, and every morphism $f$ from a perfect MV-algebra to a semisimple MV-algebra factors through either $\initial$ or $\terminal$. This naturally leads to the question of whether these two categories form a non-pointed version of a torsion theory. The appropriate notion to answer this question is that of a \emph{pretorsion theory}, which was introduced in the recent work \cite{facchini2020pretorsion}. This concept provides a generalization of the classical torsion theories, and it has been successfully applied in various areas of mathematics.\\
Additionally, we will prove that the subcategory $s\mv \subseteq \mv$ is not only reflective, but also that the adjunction $S \dashv i$ (with $S$ representing the reflector and $i$ the inclusion) is admissible for categorical Galois theory \cite{janelidze1990pure} with respect to the class of all arrows (and also with respect to the class of regular epimorphisms). By studying this Galois structure, we will be able to characterize the trivial, normal, and central extensions related to it.\\
Finally, it should be noted that the subcategory of the category of regular epimorphisms in $\mv$ whose objects are the central extensions is reflective. The corresponding Galois structure is admissible with respect to double extensions, which opens the way to studying higher-dimensional central extensions.\\
The paper is organized as follows.\\
\indent In Section \ref{Preliminaries2} we recall some classical facts about MV-algebras and we focus on the notion of ideal of an MV-algebra.\\
\indent In Section \ref{Protomodularity, Arithmeticity, and Centralizers} we show that the category of MV-algebras is a protomodular and arithmetical \cite{aritmetica} category. Since $\mv$ is a variety of universal algebras, we do this by providing explicit descriptions of the protomodularity and arithmeticity terms. Moreover, we prove that all subobjects in the category $\pt{B}{\mv}=(\mv/B)\backslash id_B$ (i.e.\ the coslice over $id_B$ of the slice of $\mv$ over $B$) have centralizers.\\
\indent In Section \ref{Preliminaries3} we will review the necessary preliminary concepts required to understand the remaining part of this work. Specifically, we will focus on the concepts of pretorsion theories, categorical Galois theory and factorization systems.\\
\indent In Section \ref{Galois Theory for MV-Algebras} we will delve into a detailed study of the adjunction determined by the reflective subcategory $s\mv$ from the perspective of categorical Galois theory. Additionally, we will describe the commutators defined by this Galois structure. Finally, we will investigate some properties of the functor $S$ and show how the Galois structure induces a stable factorization system on $\mv$.\\
\indent In Section \ref{Higher Galois Theory for MV-Algebras} we will study the higher-dimensional normal and central extensions relative to the Galois structure determined by the adjunction defined by the subcategory of regular epimorphisms whose objects are the central extensions. This in-depth analysis will allow us to define the commutator of two ideal subalgebras relatively to this Galois structure.

\section{Preliminaries}\label{Preliminaries2}

In this section, we will recall the fundamental properties of MV-algebras, starting from their definition. We will examine how new operations can be derived from those introduced in the definition. Additionally, we will recall how a partial order can be naturally defined on any MV-algebra. Finally, we will focus on the concept of ideal in an MV-algebra.
\begin{definition}

An \emph{MV-algebra} is an algebra $(A, \oplus, \lnot, 0)$ with a binary operation $\oplus$, a unary operation $\lnot$, and a constant $0$ satisfying the following equations:
\begin{itemize}
    \item[M1)] $x \oplus (y \oplus z)=(x \oplus y) \oplus z$;
    \item[M2)] $x \oplus y= y \oplus x$;
    \item[M3)] $x \oplus 0=x$;
    \item[M4)] $\lnot \lnot x=x$;
    \item[M5)] $x \oplus \lnot 0= \lnot 0$;
    \item[M6)] $\lnot (\lnot x \oplus y) \oplus y=\lnot (\lnot y \oplus x) \oplus x$.

\end{itemize}
A morphism between two MV-algebras $(A, \oplus, \lnot, 0)$ and $(B, \oplus, \lnot, 0)$ is a map $f \colon  A \rightarrow B$ satisfying the following conditions, for every $x,y \in A$:
\begin{itemize}
    \item[H1)] $f(0)=0$;
    \item[H2)] $f(x \oplus y)=f(x) \oplus f(y)$;
    \item[H3)] $f(\lnot x)= \lnot f(x)$.
\end{itemize}
The category $\mv$ is the category whose objects are the MV-algebras and whose arrows are the morphisms between them.

\end{definition}

We will often denote a model of an algebraic theory $(A,\mathcal{O})$ (where $\mathcal{O}$ is the set of internal operations) with the underlying set $A$. So, for instance, we will indicate an MV-algebra $(A, \oplus, \lnot, 0)$ simply with $A$. Moreover, we will consider the $\lnot$ operation more binding than the $\oplus$ operation.\\

Given an MV-algebra $A$, we define the constant $1$ and the binary operations $\odot$, $\ominus$, $\rightarrow$ and $d$ as follows:
\begin{itemize}
    \item $1 \coloneqq \lnot 0$;
    \item $x \odot y \coloneqq \lnot (\lnot x \oplus \lnot y)$;
    \item $x \rightarrow y \coloneqq \lnot x \oplus y$;
    \item $x \ominus y \coloneqq \lnot (\lnot x \oplus y)=\lnot (x \rightarrow y)=x \odot \lnot y$;
    \item $d(x,y) \coloneqq (x \ominus y) \oplus (y \ominus x)$ (called \emph{distance} between $x$ and $y$).
\end{itemize}

 We will consider the $\lnot$ operation more binding than the $\odot$, $\ominus$, and $\rightarrow$ operations.\\

With respect to the new operations, we get

\begin{itemize}
    \item $\lnot 1 =0$;
    \item $x \oplus y= \lnot (\lnot x \odot \lnot y)$;
    \item $x \oplus 1=1$;
    \item $(x \ominus y) \oplus y=(y \ominus x) \oplus x$;
    \item $x \rightarrow x= \lnot x \oplus x=1$;
    \item $x=y$ if and only if $d(x,y)=0$ (see  \cite{cignoli2013algebraic}, Proposition 1.2.5).
\end{itemize}

We recall now some known facts that will be useful in the next sections.

\begin{lemma}[\cite{cignoli2013algebraic}, Lemma 1.1.2]

Let $A$ be an MV-algebra and $x, y \in A$. Then the following conditions are equivalent:
\begin{itemize}
    \item $\lnot x \oplus y=1$;
    \item $x \odot \lnot y=0$;
    \item $y=(y \ominus x) \oplus x$;
    \item there exists an element $z \in A$ such that $y=x \oplus z$.
\end{itemize}

\end{lemma}

 Given an MV-algebra $A$ and two elements $x,y \in A$, we say that $$x \leq y$$ if and only if $x$ and $y$ satisfy one of the above equivalent conditions. It can be shown that $\leq$ defines a partial order on $A$ (a proof of this fact can be found in \cite{cignoli2013algebraic}). It is clear that every morphism of MV-algebras preserves the partial order defined above. 

\begin{lemma}[\cite{cignoli2013algebraic}, Lemma 1.1.4]

In every MV-algebra A the natural order $\leq$ has the following properties:

\begin{itemize}
    \item $x \leq y$ if and only if $\lnot y \leq \lnot x$;
    \item if $x \leq y$ then, for all $z \in A$, $x \oplus z \leq y \oplus z$ and $x \odot z \leq y \odot z$;
    \item $x \odot y \leq z$ if and only if $x \leq \lnot y \oplus z$.
\end{itemize}

\end{lemma}

\begin{proposition}[\cite{cignoli2013algebraic}, Proposition 1.1.5, 1.1.6, and 1.5.1]\label{inf}
On every MV-algebra $A$ the natural order determines a lattice structure. Specifically, for every $x, y \in A$, the join $x \lor y$ and the meet $x \land y$ are given by $$x \lor y \coloneqq (x \ominus y) \oplus y \textit{ and } x \land y \coloneqq x \odot (\lnot x \oplus y).$$
Moreover, the underlying lattice structure on $A$ is distributive.\\
Finally, the following equations hold:
\begin{itemize}
    \item $x \odot (y \lor z)=(x \odot y) \lor (x \odot z)$;
    \item $x \oplus (y \land z)=(x \oplus y) \land (x \oplus z)$;
    \item $x \odot (y \land z)=(x \odot y) \land (x \odot z)$;
    \item $x \oplus (y \lor z)=(x \oplus y) \lor (x \oplus z)$.
\end{itemize}
%%per le due citazione migliore
\end{proposition}

Given a morphism of MV-algebras $h \colon A \rightarrow B$, one can define the \emph{kernel} of $h$ as $$\kerr(h) \coloneqq \{ x \in A \, | \, h(x)=0 \}.$$ Kernels of morphisms can be characterized as \emph{ideals}. A subset $I \subseteq A$ of an MV-algebra $A$ is an ideal if and only if for every $x,y \in A$:
\begin{itemize}
    \item $0 \in I$;
    \item $x,y \in I$ implies $x \oplus y \in I$;
    \item $x \in I$ and $y \leq x$ implies $y \in I$.
\end{itemize}
An ideal $I$ is said to be \emph{proper} if $I \neq A$.\\
The dual notion is the one of \emph{filter}. A subset $F \subseteq A$ of an MV-algebra $A$ is a filter if and only if for every $x,y \in A$:
\begin{itemize}
    \item $1 \in F$;
    \item $x,y \in F$ implies $x \odot y \in F$;
    \item $x \in F$ and $y \geq x$ implies $y \in F$.
\end{itemize}
A filter $F$ is said to be \emph{proper} if $F \neq A$.
One can easily prove that $F$ is a filter of $A$ if and only if $\lnot F \coloneqq \{ a \in A \, | \, \lnot a \in F \}$ is an ideal of $A$.\\ 
It is known that every kernel is an ideal (we recall that every morphism preserves the order). Given an ideal $I \subseteq A$ we can define a congruence $\sim_I$ on $A$ in the following way: for every $x,y \in A$, $x \sim_I y$ if and only if $d(x,y) \in I$. One can show that the kernel of the quotient projection $\pi \colon A \rightarrow A/\sim_I$ is exactly $I$. This procedure establishes a one-to-one correspondence between kernels of morphisms with domain $A$ and ideals of $A$. Additionally, it is not difficult to see that a morphism of MV-algebras $h$ is injective if and only if $\kerr(h)=\{ 0 \}$.\\
Finally, given an MV-algebra $A$ and a non-empty subset $S \subseteq A$, the ideal generated by $S$ (i.e.\ the smallest ideal containing $S$) exists and it is $$\langle S \rangle =\{ x \in A\, | \, \exists s_1, \dots, s_n \in S \text{ such that } x \leq s_1 \oplus \dots \oplus s_n \}.$$

\begin{lemma}\label{ideali}

Let $f \colon A \rightarrow B$ be a morphism of MV-algebras and $I \subseteq A$ an ideal. The restriction of $f$ to $I$ is injective if and only if $\kerr(f) \cap I = \{ 0 \}$.

\begin{proof}

If the restriction of $f$ is injective and we take an element $a \in \kerr(f) \cap I$ we get $f(a)=0=f(0)$ and, since $0 \in I$, we obtain $a=0$. Conversely, if $\kerr(f) \cap I = \{ 0 \}$ and we consider two elements $a,b \in I$ such that $f(a)=f(b)$ we get that $f(d(a,b))=d(f(a),f(b))=0$ and so $d(a,b) \in \kerr(f) \cap I = \{ 0 \}$ i.e.\ $a=b$; thus the restriction of $f$ is injective.
\end{proof}

\end{lemma}

\section{Protomodularity, Arithmeticity, and Centralizers} \label{Protomodularity, Arithmeticity, and Centralizers}
In this section, we will study some categorical-algebraic properties of $\mv$. Since $\mv$ is a variety of universal algebras, we know that it is a Barr-exact category \cite{barr}. It is known that $\mv$ is also arithmetical \cite{aritmetica} and protomodular \cite{bourn1991normalization}; for proofs of these facts, see, respectively, \cite{idziak1984lattice} and  \cite{lapenta2022relative}. In this section, we will offer alternative proofs of these results by exhibiting explicit terms of arithmeticity and protomodularity. Given a category $\cat{C}$ and an object $B$ of it, $\pt{\cat{C}}{B}$ denotes the category whose objects are the split epimorphisms over $B$ with a fixed splitting and the morphisms are the commutative triangles between these data (or, as previously mentioned in the introduction, $\pt{\cat{C}}{B}$ is the coslice over $id_B$ of the slice of $\cat{C}$ over $B$). A finitely complete category $\cat{C}$ is \emph{protomodular} if, for every arrow $f \colon A \rightarrow B$, the change-of-base functor $f^* \colon \pt{\cat{C}}{B} \rightarrow \pt{\cat{C}}{A}$ is conservative; this property, in the pointed case, is equivalent to requiring the validity of the split short five lemma. We will then prove how protomodularity allows us to identify the conditions under which certain commutative squares in $\mv$ are pullbacks. A Barr-exact category with coequalizers $\mathbb{C}$ is \emph{arithmetical} if it is a Mal'tsev \cite{carboni1991diagram} category (i.e.\ finitely complete and such that every internal reflexive relation is an internal equivalence relation) and, for any object $B$ of $\mathbb{C}$, the lattice $\Eq(B)$ of internal equivalence relations on $B$ is distributive. Finally, we will provide a more detailed analysis of the properties of the categories $\pt{\mv}{B}$, showing in particular that in these categories every subobject has a centralizer.\\
\begin{comment}
Specifically, we will give an explicit description of the terms of protomodularity and arithmeticity for $\mv$.
\end{comment}

To determine whether a variety $\mathbb{V}$ is protomodular, we can use Theorem 1.1 of \cite{varprot}. This theorem states that $\mathbb{V}$ is protomodular if and only if it has $0$-ary terms $e_1, \dots, e_n$, binary terms $t_1, \dots, t_n$, and an $(n+1)$-ary term $t$ satisfying the identities $$t(x,t_1(x,y), \dots ,t_n(x,y))=y \text{ and } t_i(x,x)=e_i$$ for all $i=1, \dots ,n$.\\

In the recent work \cite{lapenta2022relative}, the authors make use of the results obtained for the variety of hoops to prove, among other things, that the variety of MV-algebras is protomodular. Here, we present an alternative proof of this fact by exhibiting different protomodularity terms compared to those introduced in the aforementioned work.

\begin{proposition}

$\mv$ is a protomodular category. 

\begin{proof}

We define $t_1(x,y) \coloneqq x \ominus y$, $t_2(x,y) \coloneqq x \oplus \neg y$, and $t(x,y,z) \coloneqq x \oplus (y \odot z)$. Clearly, one has $$t_1(x,x)=x \ominus x=0 \text{ and } t_2(x,x)=x \oplus \neg x=1.$$ Moreover, the following equality holds: $$t(t_1(x,y),t_2(x,y),y)=(x \ominus y) \oplus ((x \oplus \neg y) \odot y)=x;$$
a proof of the last equality can be found in \cite{cignoli2013algebraic}, Proposition 1.6.2.
\end{proof}

\end{proposition}

In their work \cite{lapenta2022relative}, the authors show that $\mv$ is a protomodular category by constructing two binary terms, namely $\alpha_1(x,y)$ and $\alpha_2(x,y)$, as well as a ternary term, $\theta(x,y,z)$. However, we observe that, in their case, the equalities $\alpha_1(x,x)=1$ and $\alpha_2(x,x)=1$ hold. Consequently, it is clear that our own protomodularity terms differ from theirs, as we mentioned before.

\begin{lemma}\label{Lemma dei pullback}

Consider a commutative square in $\mv$
% https://q.uiver.app/?q=WzAsNCxbMCwwLCJBIl0sWzEsMCwiQiJdLFswLDEsIkMiXSxbMSwxLCJEIl0sWzAsMSwiZiIsMCx7InN0eWxlIjp7ImhlYWQiOnsibmFtZSI6ImVwaSJ9fX1dLFsxLDMsImsiLDAseyJzdHlsZSI6eyJoZWFkIjp7Im5hbWUiOiJlcGkifX19XSxbMiwzLCJnIiwyLHsic3R5bGUiOnsiaGVhZCI6eyJuYW1lIjoiZXBpIn19fV0sWzAsMiwiaCIsMix7InN0eWxlIjp7ImhlYWQiOnsibmFtZSI6ImVwaSJ9fX1dXQ==
\[\begin{tikzcd}
	A & B \\
	C & D,
	\arrow["f", two heads, from=1-1, to=1-2]
	\arrow["k", from=1-2, to=2-2]
	\arrow["g"', two heads, from=2-1, to=2-2]
	\arrow["h"', from=1-1, to=2-1]
\end{tikzcd}\]
where the horizontal arrows are regular epimorphisms. Hence, denoting with $\langle h, f \rangle$ the unique arrow induced by the universal property of the pullback
% https://q.uiver.app/?q=WzAsNSxbMCwwLCJBIl0sWzIsMSwiQiJdLFsxLDIsIkMiXSxbMiwyLCJEIl0sWzEsMSwiUCJdLFswLDEsImYiLDAseyJzdHlsZSI6eyJoZWFkIjp7Im5hbWUiOiJlcGkifX19XSxbMSwzLCJrIiwwLHsic3R5bGUiOnsiaGVhZCI6eyJuYW1lIjoiZXBpIn19fV0sWzIsMywiZyIsMix7InN0eWxlIjp7ImhlYWQiOnsibmFtZSI6ImVwaSJ9fX1dLFswLDIsImgiLDIseyJzdHlsZSI6eyJoZWFkIjp7Im5hbWUiOiJlcGkifX19XSxbNCwyLCJcXHBpX0MiLDIseyJzdHlsZSI6eyJoZWFkIjp7Im5hbWUiOiJlcGkifX19XSxbNCwxLCJcXHBpX0IiXSxbNCwzLCIiLDEseyJzdHlsZSI6eyJuYW1lIjoiY29ybmVyIn19XSxbMCw0LCJcXGxhbmdsZSBoLGYgXFxyYW5nbGUiLDJdXQ==
\[\begin{tikzcd}
	A \\
	& P & B \\
	& C & D,
	\arrow["f", bend left, two heads, from=1-1, to=2-3]
	\arrow["k", from=2-3, to=3-3]
	\arrow["g"', two heads, from=3-2, to=3-3]
	\arrow["h"', bend right, from=1-1, to=3-2]
	\arrow["{\pi_C}"', from=2-2, to=3-2]
	\arrow["{\pi_B}", two heads, from=2-2, to=2-3]
	\arrow["\lrcorner"{anchor=center, pos=0.125}, draw=none, from=2-2, to=3-3]
	\arrow["{\langle h,f \rangle}", from=1-1, to=2-2]
\end{tikzcd}\]
we get:
\begin{itemize}
    \item[i)] $\langle h, f \rangle $ is injective if and only the restriction of $h$ (considered as map) $h \colon \kerr(f) \rightarrow \kerr(g)$ is injective;
    \item[ii)] $\langle h, f \rangle $ is surjective if and only the restriction of $h$ (considered as map) $h \colon \kerr(f) \rightarrow \kerr(g)$ is surjective.
\end{itemize}
\begin{proof}
$i)$ $(\Rightarrow)$ By assumption $\langle h,f \rangle$ is injective, therefore $\{0 \} = \kerr(\langle h,f \rangle) = \kerr(f) \cap \kerr(h)$ and so, thanks to Lemma \ref{ideali}, we conclude that the restriction of $h$ is injective.\\
$(\Leftarrow)$ Consider an element $a \in A$ such that $\langle h,f \rangle (a)=(0,0)$; then $ a \in \kerr(f) \cap \kerr(h)$ and so $a=0$, since the restriction of $h$ is injective.\\
$ii)$ $(\Rightarrow)$ Fix an element $c \in \kerr(g)$; we know that $(c,0) \in P$ since $k(0)=0=g(c)$; hence, there exists an element $a \in A$ such that $\langle h,f \rangle (a)=(c,0)$ ($\langle h,f \rangle$ is surjective). Therefore, $a \in \kerr(f)$ and $h(a)=c$, i.e.\ the restriction of $h$ is surjective.\\
($\Leftarrow$) Given an element $(c,b) \in P$ (i.e.\ $c \in C$, $b \in B$, and $g(c)=k(b)$) there exists an element $a \in A$ such that $f(a)=b$ ($f$ is surjective). Now, we observe that $gh(a)=kf(a)=k(b)=g(c)$, and so we deduce $c \ominus h(a) \in \kerr(g)$ and $h(a) \ominus c \in \kerr(g)$. Therefore, since the restriction of $h$ is surjective, there exist $a_1,a_2 \in \kerr(f)$ such that $h(a_1)=c \ominus h(a) \in \kerr(g)$ and $h(a_2)=h(a) \ominus c \in \kerr(g)$ (for our aim, it is most useful to keep in mind that $h(\lnot a_2)=c \oplus \lnot h(a)$). We recall, from \cite{cignoli2013algebraic}, Proposition 1.6.2, that the following equality holds: $$(c \ominus h(a)) \oplus ((c \oplus \lnot h(a)) \odot h(a))=c.$$ This equality can be reformulated as $$h(a_1) \oplus (h(\lnot a_2) \odot h(a))=c,$$ so $h(a_1 \oplus (\lnot a_2 \odot a))=c$ and $f(a_1 \oplus (\lnot a_2 \odot a))=f(a_1) \oplus (\lnot f(a_2) \odot f(a))=0 \oplus (1 \odot f(a))=b$; we have proved that $\langle h, f \rangle$ is surjective.
\end{proof}

\end{lemma}

The previous lemma yields the following direct consequence:

\begin{corollary}\label{Lemma dei pullback 2}

A commutative square in $\mv$
% https://q.uiver.app/?q=WzAsNCxbMCwwLCJBIl0sWzEsMCwiQiJdLFswLDEsIkMiXSxbMSwxLCJEIl0sWzAsMSwiZiIsMCx7InN0eWxlIjp7ImhlYWQiOnsibmFtZSI6ImVwaSJ9fX1dLFsxLDMsImsiLDAseyJzdHlsZSI6eyJoZWFkIjp7Im5hbWUiOiJlcGkifX19XSxbMiwzLCJnIiwyLHsic3R5bGUiOnsiaGVhZCI6eyJuYW1lIjoiZXBpIn19fV0sWzAsMiwiaCIsMix7InN0eWxlIjp7ImhlYWQiOnsibmFtZSI6ImVwaSJ9fX1dXQ==
\[\begin{tikzcd}
	A & B \\
	C & D,
	\arrow["f", two heads, from=1-1, to=1-2]
	\arrow["k", from=1-2, to=2-2]
	\arrow["g"', two heads, from=2-1, to=2-2]
	\arrow["h"', from=1-1, to=2-1]
\end{tikzcd}\]
where the horizontal arrows are regular epimorphism, is a pullback if and only if the restriction of $h$ (considered as map) $h \colon \kerr(f) \rightarrow \kerr(g)$ is bijective.

\end{corollary}

It is well-known that in a pointed protomodular category, a commutative square in which two parallel arrows are regular epimorphisms forms a pullback if and only if the kernels of the two regular epimorphisms are isomorphic. Thus, the statement in Corollary \ref{Lemma dei pullback 2} could emerge as a reformulation of this result in the category of MV-algebras, which is not pointed. Importantly, it should be noticed that Corollary \ref{Lemma dei pullback 2} does not directly stem from the known results in the pointed case. Specifically, in $\mv$, for an arbitrary arrow $f$, $\kerr(f)$ is not, in general, a subobject of the domain of $f$.\\

In \cite{pixley1963distributivity}, Theorem 2, the author proved that a variety $\cat{V}$ is arithmetical if and only if there exists a ternary term $r(x,y,z)$ such that $$r(x,x,z)=z \text{, } r(x,y,y)=x \text{, and } r(x,y,x)=x$$
for every object $X$ and for every $x,y,z \in X$. The fact that the category $\mv$ is arithmetical can be seen as a consequence of a result by Idziak in \cite{idziak1984lattice}. We provide a different proof of this by exhibiting an explicit term of arithmeticity:

\begin{proposition}

$\mv$ is an arithmetical category.

\begin{proof}

We define $$p(x,y,z) \coloneqq ((x \rightarrow y) \rightarrow z) \land ((z \rightarrow y) \rightarrow x)$$ and $$t(x,y,z) \coloneqq (y \rightarrow (x \land z)) \land (x \lor z).$$ We observe that 
\begin{align*}
p(x, x, z)&=z \land((z \rightarrow x) \rightarrow x)=z \land(\neg(\neg z \oplus x) \oplus x)\\
&=z \land(\neg(\neg x \oplus z) \oplus z)=z,\\
p(x, y, y) &=((x \rightarrow y) \rightarrow y) \land x=(\neg(\neg x \oplus y) \oplus y) \land x\\
&=(\neg(\neg y \oplus x) \oplus x) \land x=x, \text{ and }\\
p(x, y, x)&=((x \rightarrow y) \rightarrow x) \land((x \rightarrow y) \rightarrow x)\\
&=((x \rightarrow y) \rightarrow x)=\neg(\neg x \oplus y) \oplus x \geq x.    
\end{align*}
Moreover, we have
\begin{align*}
t(x,x,z)&=(x \rightarrow (x \land z)) \land (x \lor z)=((x \rightarrow x) \land (x \rightarrow z)) \land (x \lor z)\\
&=(x \rightarrow z) \land (x \lor z) \geq z,\\
t(x,y,y)&=(y \rightarrow (x \land y)) \land (x \lor y)=(y \rightarrow x) \land (x \lor y) \geq x, \text{ and} \\
t(x,y,x)&=(y \rightarrow x) \land x=x.  
\end{align*}
Therefore, the term $$r(x,y,z) \coloneqq p(x,y,z) \land t(x,y,z)$$ satisfies $$r(x,x,z)=z \text{, } r(x,y,y)=x \text{, and } r(x,y,x)=x.\eqno{\qedhere}$$
\end{proof}

\end{proposition}
In the final part of this section we study, from a categorical point of view, the commutativity of subobjects in the category $\pt{B}{\mv}$, for every MV-algebra $B$.\\
To introduce the topic, we mention some known results related to the category of groups. Given a group $G$ and two subgroups $A,B \leq G$, the condition that, for every $a \in A$ and $b \in B$, $ab=ba$ can be reformulated in this equivalent way: there exists a group homomorphism $\varphi \colon  A \times B \rightarrow G$ making the following diagram commutative:
% https://q.uiver.app/?q=WzAsNCxbMCwwLCJBIl0sWzIsMCwiQiJdLFsxLDEsIkEgXFx0aW1lcyBCIl0sWzEsMywiRyJdLFswLDIsIihpZF9BLDApIiwxXSxbMSwyLCIoMCxpZF9CKSIsMV0sWzAsMywiIiwwLHsic3R5bGUiOnsidGFpbCI6eyJuYW1lIjoiaG9vayIsInNpZGUiOiJ0b3AifX19XSxbMSwzLCIiLDIseyJzdHlsZSI6eyJ0YWlsIjp7Im5hbWUiOiJob29rIiwic2lkZSI6ImJvdHRvbSJ9fX1dLFsyLDMsIlxccGhpIiwxXV0=
\[\begin{tikzcd}
	A && B \\
	& {A \times B} \\
	\\
	& G.
	\arrow["{(id_A,0)}"{description}, from=1-1, to=2-2]
	\arrow["{(0,id_B)}"{description}, from=1-3, to=2-2]
	\arrow[bend right, hook, from=1-1, to=4-2]
	\arrow[bend left, hook', from=1-3, to=4-2]
	\arrow["\varphi", from=2-2, to=4-2]
\end{tikzcd}\]
Moreover, observing that $\varphi(a,b)=\varphi(a,1)\varphi(1,b)=ab$, we conclude that $\varphi$ must be the group product and, therefore, it is necessarily unique. Hence, it makes sense to place ourselves in a context in which a morphism  $\varphi$ of this type is unique. This reasoning justifies the following definition:

\begin{definition}[\cite{unitalcat}] A pointed category $\mathbb{C}$ with finite products is \emph{unital} if for all objects $X,Y$ the pair of morphisms $(id_X,0) \colon X \rightarrow X \times Y$, $(0,id_Y) \colon Y \rightarrow X \times Y$ is jointly extremally epimorphic.

\end{definition}

To be more explicit, a pair of arrows $f \colon  A \rightarrow B$ and $g \colon  C \rightarrow B$ of a category $\mathbb{C}$ is said to be jointly extremally epimorphic when for every commutative diagram

% https://q.uiver.app/?q=WzAsNCxbMCwxLCJBIl0sWzEsMSwiQiJdLFsyLDEsIkMiXSxbMSwwLCJNIl0sWzAsMSwiZiIsMl0sWzIsMSwiZyJdLFswLDMsImYnIl0sWzIsMywiZyciLDJdLFszLDEsIm0iXV0=
\[\begin{tikzcd}
	& M \\
	A & B & C
	\arrow["f"', from=2-1, to=2-2]
	\arrow["g", from=2-3, to=2-2]
	\arrow["{f'}", from=2-1, to=1-2]
	\arrow["{g'}"', from=2-3, to=1-2]
	\arrow["m", from=1-2, to=2-2]
\end{tikzcd}\]

if $m$ is a monomorphism, then $m$ is an isomorphism.\\

It has been shown in \cite{semi} that every semi-abelian category is unital.\\

We are ready to mention the generalized notion of commutativity between subobjects.

\begin{definition}[\cite{huq}]
Let $\cat{C}$ be a unital category. Two subobjects $a \colon  A \mono X$ and $b \colon  B \mono X$ of $X$ are said to \emph{cooperate} (or \emph{commute in the sense of Huq}, and we write $[a,b]=0$) if there exists a (necessarily unique) morphism $\varphi \colon  A \times B \rightarrow X$ (called \emph{cooperator}) such that the following diagram commutes:
% https://q.uiver.app/?q=WzAsNCxbMCwwLCJBIl0sWzEsMywiWCJdLFsyLDAsIkIiXSxbMSwxLCJBIFxcdGltZXMgQiJdLFswLDEsIiIsMCx7ImN1cnZlIjoxLCJzdHlsZSI6eyJ0YWlsIjp7Im5hbWUiOiJob29rIiwic2lkZSI6InRvcCJ9fX1dLFsyLDEsIiIsMix7ImN1cnZlIjotMSwic3R5bGUiOnsidGFpbCI6eyJuYW1lIjoiaG9vayIsInNpZGUiOiJib3R0b20ifX19XSxbMCwzLCJpX0EiLDFdLFsyLDMsImlfQiIsMV0sWzMsMSwiXFxleGlzdCAhIFxccGhpIiwxXV0=
\[\begin{tikzcd}
	A && B \\
	& {A \times B} \\
	\\
	& X.
	\arrow[ bend right, tail , "{a}"', from=1-1, to=4-2]
	\arrow[ bend left, tail , "{b}", from=1-3, to=4-2]
	\arrow["{(id_A,0)}"{description}, from=1-1, to=2-2]
	\arrow["{(0,id_B)}"{description}, from=1-3, to=2-2]
	\arrow["{\exists  \varphi}", from=2-2, to=4-2]
\end{tikzcd}\]
Given a subobject $a \colon  A \mono X$, the \emph{centralizer} of $a$ in $X$, if it exists, is the greatest subobject of $X$ that cooperates with $a$.
\end{definition}

To show that in $\pt{B}{\mv}$ there are centralizers of subobjects, we need to recall the notion of lattice-ordered abelian group and its link with MV-algebras. A \emph{lattice-ordered abelian group} is an algebraic structure of signature $\{ +, 0, -, \lor, \land \}$ satisfying the axioms of abelian groups, the axioms of lattices, and the axioms related to the distributivity of the group operation over both the lattice operations: $$x+(y \lor z)=(x+y) \lor (x+z) \text{ and } x+(y \land z)=(x+y) \land (x+z).$$ We denote with $\lab$ the category whose objects are lattice-ordered abelian groups and whose arrows are maps between lattice-ordered abelian groups which preserves the operations. Given a lattice-ordered abelian group $G$ we can define, for every element $x$ of $G$, $|x| \coloneqq x \lor -x$. An \emph{order-unit} $u$ of $G$ is an element $0 \leq u \in G$ satisfying the following property: for every $x \in G$, there exists a natural number $n \in \mathbb{N}$ such that $|x| \leq nu$. We denote with $\ulab$ the category whose objects are the pairs $(G,u)$, where $G$ is a lattice-ordered abelian group and $u$ is an order-unit of $G$, and whose arrows are the maps which preserve the operations and the distinguished order-unit. Given an object $(G,u)$ of $\ulab$ we consider $$[0,u] \coloneqq \{ x \in G \, | \, 0 \leq x \leq u \};$$on $[0,u]$ the following operations are defined: $x \oplus y \coloneqq (x+y) \land u$ and $\lnot x \coloneqq u-x$. The structure $([0,u], \oplus, \lnot, 0)$ is an MV-algebra, denoted with $\Gamma(G,u)$; moreover, for every arrow $h \colon (G,u) \rightarrow (H,v)$ in $\ulab$, the restriction of $h$ to $[0,u]$ (denoted with $\Gamma(h)$) is a morphism of MV-algebras between $[0,u]$ and $[0,v]$. 

\begin{theorem}[\cite{mundici1986interpretation}, Theorem 3.9]

The assignment defined by $\Gamma$ establishes an equivalence of categories between $\ulab$ and $\mv$.

\end{theorem}

We know that $\ulab$ is complete and cocomplete, since it is equivalent to a variety of universal algebras. We want to describe finite limits in $\ulab$. We prove that they are computed as in $\lab$. We start dealing with equalizers. Let us consider a diagram in $\ulab$ of the form 
% https://q.uiver.app/?q=WzAsMixbMCwwLCIoWCx1KSJdLFsxLDAsIihZLHYpIl0sWzAsMSwiZiIsMCx7Im9mZnNldCI6LTF9XSxbMCwxLCJnIiwyLHsib2Zmc2V0IjoxfV1d
\[\begin{tikzcd}
	{(X,u)} & {(Y,v);}
	\arrow["f", shift left=1, from=1-1, to=1-2]
	\arrow["g"', shift right=1, from=1-1, to=1-2]
\end{tikzcd}\]
We define $E \coloneqq \{ x \in X \, | \, f(x)=g(x) \}$. We know that $E$ inherits the lattice-ordered abelian group operations from $X$ (limits in $\lab$ are computed as in $\mathbb{S}et$ since $\lab$ is a variety). Moreover, since $f(u)=v=g(u)$, we have $u \in E$. So, given an arrow $k \colon (H,h) \rightarrow (X,u)$ of $\ulab$ such that $fk=gk$, then, since $E$ is the equalizer of $f$ and $g$ in $\lab$, $k$ factors through the inclusion of $E$ in $X$; moreover, $k(h)=u$ and so we can conclude that $(E,u) \hookrightarrow (X,u)$ is the equalizer of $f$ and $g$ in $\ulab$. Let us take a look at the products now. We consider two objects $(X,u)$ and $(Y,v)$ of $\ulab$. We prove that $(X \times Y, (u,v))$ is an object of $\ulab$ (where the operations on $X \times Y$ are defined component-wise); in other terms, we have to show that $(u,v)$ is an order-unit. So, fix an element $(x,y) \in X \times Y$. Then, there exist $n_1, n_2 \in \mathbb{N}$ such that $  |  x  |   \leq n_1u$ and $  |  y  |   \leq n_2v$. Thus, taking $n$ as the maximum between $n_1$ and $n_2$, we get $  |  (x,y)  |  =(  |  x  |  ,  |  y  |  ) \leq (nu,nv)=n(u,v)$. Hence, applying similar reasoning to the one seen for equalizers, we obtain that the products in $u\lab$ are computed as in $\lab$.

\begin{proposition}

$\pt{(B,v)}{\ulab}$ is a unital category.

\begin{proof}

Since $\ulab$ is arithmetical it is also Mal'tsev and so we can apply Proposition 10 of \cite{unitalcat} to get that $\pt{(B,v)}{\ulab}$ is unital for every object $(B,v)$.
\end{proof}

\end{proposition}

\begin{proposition}

Let $(B,v)$ be an object of $\ulab$. In the category $\pt{(B,v)}{\ulab}$ subobjects have centralizers.

\begin{proof}

Let $((A,u),p,s)$ be an object of $\pt{(B,v)}{\ulab}$ and let us suppose, without loss of generality, that $s \colon (B,v) \rightarrow (A,u)$ is the inclusion (in particular, this implies $u=v$). We observe that $(A,p,s)$ is an object of $\pt{B}{\lab}$. We define $K \coloneqq \{ k \in A \, | \, p(k)=0 \}$ and we observe that $K$ is a subalgebra of $A$ in $\lab$. Hence, applying the results from \cite{CappellettiLgruppi}, Proposition 7.3, we obtain that $A$ is isomorphic as a lattice-ordered abelian group to the semi-direct product $K \rtimes B$, whose operations are defined by $$(k_1,b_1)+(k_2,b_2)=(k_1+k_2,b_1+b_2)$$
and $$(k_1,b_1)\lor(k_2,b_2)=(((k_1+b_1) \lor (k_2+b_2))-(b_1 \lor b_2),b_1 \lor b_2).$$ In other words, the object $(A,p,s)$ is isomorphic to
% https://q.uiver.app/?q=WzAsMixbMCwwLCJLIFxccnRpbWVzIEIiXSxbMSwwLCJCIl0sWzAsMSwicF9CIiwyLHsib2Zmc2V0IjoxfV0sWzEsMCwiaV9CIiwyLHsib2Zmc2V0IjoxfV1d
\[\begin{tikzcd}
	{K \rtimes B} & B
	\arrow["{p_B}"', shift right=1, from=1-1, to=1-2]
	\arrow["{i_B}"', shift right=1, from=1-2, to=1-1]
\end{tikzcd}\]
where $p_B(k,b)=b$ and $i_B(b)=(0,b)$. The isomorphism is given by the arrow of $\ulab$ $\varphi \colon K \rtimes B \rightarrow A$ defined by $\varphi(x,b) \coloneqq x+b$. Clearly $\varphi$ induces an isomorphism in $\ulab$ between $(A,u)$ and $(K \rtimes B, (0,u))$. Hence, we can apply the same argument of Proposition 7.3 in \cite{CappellettiLgruppi} to establish the validity of the statement.
\end{proof}

\end{proposition}

\section{Pretorsion Theories, Galois Theory and Factorization Systems}\label{Preliminaries3}
The aim of this section is to provide an introduction to the foundational concepts required for studying pretorsion theories in general categories. Given a pair of full replete subcategories $(\pre{T}, \pre{F})$ of a category $\cat{C}$, the authors of \cite{facchini2020pretorsion} start by defining a new full subcategory $\zeros \coloneqq \pre{T} \cap \pre{F}$ of \emph{trivial} objects in the category $\cat{C}$. A morphism is considered to be $\zeros$-trivial if it factors through an object of $\zeros$. We denote the collection of $\zeros$-trivial (or simply trivial) morphisms as $\zeroarrows$. This allows to define $\zeros$-prekernels, $\zeros$-precokernels, and short $\zeros$-pre-exact sequences (or simply prekernels, precokernels, and pre-exact sequences).
\begin{definition}[\cite{facchini2020pretorsion}]
Let $f \colon A \rightarrow B$ be a morphism in $\cat{C}$. We say that a morphism $k \colon K \rightarrow A$ in $\cat{C}$ is
a \emph{$\zeros$-prekernel} of $f$ if the following properties are satisfied:
\begin{itemize}
    \item $fk$ is a morphism of $\zeroarrows$;
    \item whenever $e: E \rightarrow A$ is a morphism in $\cat{C}$ and $fe$ is in $\zeroarrows$, there exists a unique morphism $\varphi \colon E \rightarrow K$ in $\cat{C}$ such that $k \varphi = e$.
\end{itemize}
\end{definition}
Dually, we have the definition of $\zeros$-precokernel.\\
In \cite{facchini2020pretorsion}, the authors show that some properties known for kernels also hold for $\zeros$-prekernels. In particular, they prove that every $\zeros$-prekernel is a monomorphism, and that, for every morphism $f \colon A \rightarrow B$ in $\cat{C}$, the $\zeros$-prekernel of $f$ is unique up to a unique isomorphism. This means that if $k \colon K \rightarrow A$ and $k' \colon K' \rightarrow A$ are $\zeros$-prekernels of the same arrow $f$, then there exists a unique isomorphism $\varphi \colon K' \rightarrow K$ such that $k \varphi = k'$.\\
Clearly, the duals of the previous observations hold for $\zeros$-precokernels.
\begin{definition}[\cite{facchini2020pretorsion}]
    Let $f \colon A \rightarrow B$ and $g \colon B \rightarrow C$ be morphisms in $\cat{C}$. We say that % https://q.uiver.app/?q=WzAsMyxbMCwwLCJBIl0sWzEsMCwiQiJdLFsyLDAsIkMiXSxbMCwxLCJmIl0sWzEsMiwiZyJdXQ==
\[\begin{tikzcd}
	A & B & C
	\arrow["f", from=1-1, to=1-2]
	\arrow["g", from=1-2, to=1-3]
\end{tikzcd}\]
is a \emph{short $\zeros$-pre-exact sequence} in $\cat{C}$ if $f$ is a $\zeros$-prekernel of $g$ and $g$ is a $\zeros$-precokernel of $f$.
\end{definition}
A pretorsion theory \cite{facchini2021pretorsion} in a category $\cat{C}$ is defined as a pair $\torsione$ of full and replete subcategories $\pre{T}$ and $\pre{F}$ of a category $\cat{C}$, which satisfy certain conditions. Specifically, every morphism from an object in $\pre{T}$ to an object in $\pre{F}$ must be trivial, and, for every object $A$ in $\cat{C}$, there must exist a pre-exact sequence with a torsion object in $\pre{T}$ as its left endpoint and a torsion-free object in $\pre{F}$ as its right endpoint.
This broader view allows for more flexibility in the choice of the category $\cat{C}$ and the subcategories $\pre{T}$ and $\pre{F}$. The concept of pretorsion theory can be seen as a generalization of the notion of torsion theory. In fact, when the category $\cat{C}$ is pointed, every pretorsion theory such that $\pre{T} \cap \pre{F}$ reduces to the zero object is, actually, a torsion theory.
\begin{definition}[\cite{facchini2021pretorsion}]\label{definizione preesatta}
Let $\cat{C}$ be an arbitrary category. A \emph{pretorsion theory} $(\pre{T}, \pre{F})$ in $\cat{C}$ consists of two full, replete subcategories $\pre{T}, \pre{F}$ of $\cat{C}$ satisfying the following two conditions. Set $\zeros = \pre{T} \cap \pre{F}$:
\begin{itemize}
    \item $\cat{C}(T, F) \subseteq \zeroarrows$ for every object $T \in \pre{T}$ , $F \in \pre{F}$;
    \item for every object $A$ of $\cat{C}$ there is a short $\zeros$-pre-exact sequence
    % https://q.uiver.app/?q=WzAsMyxbMCwwLCJUKEEpIl0sWzEsMCwiQSJdLFsyLDAsIkYoQSkiXSxbMCwxLCJcXHZhcmVwc2lsb25fQSJdLFsxLDIsIlxcZXRhX0EiXV0=
\[\begin{tikzcd}
	{T(A)} & A & {F(A)}
	\arrow["{\varepsilon_A}", from=1-1, to=1-2]
	\arrow["{\eta_A}", from=1-2, to=1-3]
\end{tikzcd}\]
with $T(A) \in \pre{T}$ and $F(A) \in \pre{F}$.
\end{itemize}
\end{definition}
Such a $\zeros$-pre-exact sequence is unique up to isomorphisms.\\
In \cite{facchini2021pretorsion}, the authors prove that, by fixing for each object $A$ a $\zeros$-pre-exact sequence as in Definition \ref{definizione preesatta}, every pretorsion theory gives rise to a pair of functors:
% https://q.uiver.app/?q=WzAsMTIsWzAsMCwiXFxjYXR7Q30iXSxbMiwwLCJcXHByZXtGfSJdLFswLDEsIkEiXSxbMiwxLCJGKEEpIl0sWzAsNCwiQiJdLFsyLDQsIkYoQikiXSxbNCwwLCJcXGNhdHtDfSJdLFs2LDAsIlxccHJle1R9Il0sWzQsMSwiQSJdLFs2LDEsIlQoQSkiXSxbNCw0LCJCIl0sWzYsNCwiVChCKSJdLFswLDEsIkYiXSxbMiw0LCJmIiwyXSxbMiwzLCIiLDAseyJzdHlsZSI6eyJ0YWlsIjp7Im5hbWUiOiJtYXBzIHRvIn19fV0sWzMsNSwiRihmKSJdLFs0LDUsIiIsMix7InN0eWxlIjp7InRhaWwiOnsibmFtZSI6Im1hcHMgdG8ifX19XSxbNiw3LCJUIl0sWzgsMTAsImYiLDJdLFs5LDExLCJUKGYpIl0sWzgsOSwiIiwwLHsic3R5bGUiOnsidGFpbCI6eyJuYW1lIjoibWFwcyB0byJ9fX1dLFsxMCwxMSwiIiwyLHsic3R5bGUiOnsidGFpbCI6eyJuYW1lIjoibWFwcyB0byJ9fX1dXQ==
\[\begin{tikzcd}[row sep=5pt]
	{\cat{C}} && {\pre{F}} && {\cat{C}} && {\pre{T}} \\
	A && {F(A)} && A && {T(A)} \\
	\\
	\\
	B && {F(B)} && B && {T(B),}
	\arrow["F", from=1-1, to=1-3]
	\arrow["f"', from=2-1, to=5-1]
	\arrow[maps to, from=2-1, to=2-3]
	\arrow["{F(f)}", from=2-3, to=5-3]
	\arrow[maps to, from=5-1, to=5-3]
	\arrow["T", from=1-5, to=1-7]
	\arrow["f"', from=2-5, to=5-5]
	\arrow["{T(f)}", from=2-7, to=5-7]
	\arrow[maps to, from=2-5, to=2-7]
	\arrow[maps to, from=5-5, to=5-7]
\end{tikzcd}\]
where $T(f) \colon T(A) \rightarrow T(B)$ is the unique morphism such that $f \varepsilon_A=\varepsilon_B T(f)$ and it exists since $\eta_B f \varepsilon_A \in \zeroarrows$; in a similar way, $F(f) \colon F(A) \rightarrow F(B)$ is the unique morphism such that $F(f) \eta_A = \eta_B f$ and it exists since $\eta_B f \varepsilon_A \in \zeroarrows$:
% https://q.uiver.app/?q=WzAsNixbMCwwLCJUKEEpIl0sWzEsMCwiQSJdLFsyLDAsIkYoQSkiXSxbMSwxLCJCIl0sWzIsMSwiRihCKSJdLFswLDEsIlQoQikiXSxbMSwzLCJmIl0sWzEsMiwiXFxldGFfQSJdLFszLDQsIlxcZXRhX0IiLDJdLFsyLDQsIlxcZXhpc3RzIUYoZikiXSxbMCwxLCJcXHZhcmVwc2lsb25fQSJdLFs1LDMsIlxcdmFyZXBzaWxvbl9CIiwyXSxbMCw1LCJcXGV4aXN0cyFUKGYpIiwyXV0=
\[\begin{tikzcd}
	{T(A)} & A & {F(A)} \\
	{T(B)} & B & {F(B).}
	\arrow["f", from=1-2, to=2-2]
	\arrow["{\eta_A}", from=1-2, to=1-3]
	\arrow["{\eta_B}"', from=2-2, to=2-3]
	\arrow["{\exists!F(f)}", from=1-3, to=2-3]
	\arrow["{\varepsilon_A}", from=1-1, to=1-2]
	\arrow["{\varepsilon_B}"', from=2-1, to=2-2]
	\arrow["{\exists!T(f)}"', from=1-1, to=2-1]
\end{tikzcd}\]
\begin{proposition}[\cite{facchini2021pretorsion}, Proposition 3.3]
Let $(\pre{T} ,\pre{F})$ be a pretorsion theory in a category $\cat{C}$. Then:
\begin{itemize}
    \item the functor $F \colon \cat{C} \rightarrow \pre{F}$ is a left inverse left adjoint of the inclusion functor $i_{\pre{F}} \colon \pre{F} \hookrightarrow \cat{C}$ and the unit is given by $\eta$; 
    \item the functor $T \colon \cat{C} \rightarrow \pre{T}$ is a left inverse right adjoint of the inclusion functor $i_{\pre{T}} \colon \pre{T} \hookrightarrow \cat{C}$ and the counit is given by $\varepsilon$.
\end{itemize}
\end{proposition}

In the following pages, we will review several important definitions from categorical Galois theory.

\begin{definition}

Let $\cat{C}$ be a category with pullbacks. A class $\mathcal{C}$ of arrows in $\cat{C}$ is admissible when:
\begin{itemize}
    \item every isomorphism is in $\mathcal{C}$;
    \item $\mathcal{C}$ is closed under composition;
    \item in the pullback
    % https://q.uiver.app/?q=WzAsNCxbMCwwLCJcXGJ1bGxldCJdLFsxLDAsIlxcYnVsbGV0Il0sWzEsMSwiXFxidWxsZXQiXSxbMCwxLCJcXGJ1bGxldCJdLFswLDEsImMiXSxbMywyLCJhIiwyXSxbMCwzLCJkIiwyXSxbMSwyLCJiIl0sWzAsMiwiIiwxLHsic3R5bGUiOnsibmFtZSI6ImNvcm5lciJ9fV1d
\[\begin{tikzcd}
	\bullet & \bullet \\
	\bullet & \bullet
	\arrow["c", from=1-1, to=1-2]
	\arrow["a"', from=2-1, to=2-2]
	\arrow["d"', from=1-1, to=2-1]
	\arrow["b", from=1-2, to=2-2]
	\arrow["\lrcorner"{anchor=center, pos=0.125}, draw=none, from=1-1, to=2-2]
\end{tikzcd}\]
if $a$ and $b$ are in $\mathcal{C}$ then $c$ and $d$ are in $\mathcal{C}$, too.
\end{itemize}

\end{definition}

\begin{definition}

Let $\mathcal{C}$ be an admissible class of morphisms in a category $\cat{C}$. For an object $C$ of $\cat{C}$, we write $\mathcal{C}/C$ for the following category:
\begin{itemize}
    \item the objects are the pairs $(X,f)$ where $f \colon X \rightarrow C$ in $\mathcal{C}$; 
    \item the arrows $h \colon (X,f) \rightarrow (Y,g)$ are all arrows in $\cat{C}$ such that $gh=f$.
\end{itemize}

\end{definition}

\begin{definition}[\cite{janelidze1990pure}]

A \emph{relatively admissible adjunction} consists in an adjunction $S \dashv C \colon \cat{P} \rightarrow \cat{A}$ (where $\cat{A}$ and $\cat{B}$ have pullbacks) and two admissible classes $\mathcal{A} \subseteq \cat{A}$, $\mathcal{P} \subseteq \cat{P}$ of arrows, such that
\begin{itemize}
    \item $S(\mathcal{A}) \subseteq \mathcal{P}$,
    \item $C(\mathcal{P}) \subseteq \mathcal{A}$,
    \item for every object $A$ of $\cat{A}$ the component $\eta_A$ of the unit of the adjunction $S \dashv C$ is in $\mathcal{A}$,
    \item for every object $X$ of $\cat{P}$ the component $\varepsilon_X$ of the counit of adjunction $S \dashv C$ is in $\mathcal{P}$.
\end{itemize}
We denote a relatively admissible adjunction with $(S,C,\mathcal{A},\mathcal{P})$.
\end{definition}

In order to recall the notion of an admissible adjunction, remaining consistent with the notation of previous definitions, we review the definition of the following functors: $S_A \colon \mathcal{A}/A \rightarrow \mathcal{P}/S(A)$ and $C_A \colon \mathcal{P}/S(A) \rightarrow \mathcal{A}/A$ where $A$ is an object of $\cat{A}$. One has $S_A(f \colon B \rightarrow A)  \coloneqq S(f)$ and $C_A(g \colon P \rightarrow S(A)) \coloneqq \pi_A$ where the diagram
% https://q.uiver.app/?q=WzAsNCxbMCwwLCJLIl0sWzEsMCwiQSJdLFsxLDEsIkNTKEEpIl0sWzAsMSwiQyhQKSJdLFszLDIsIkMoZykiLDJdLFsxLDIsIlxcZXRhX0EiXSxbMCwxLCJcXHBpX0EiXSxbMCwzLCJcXHBpX3tDKFApfSIsMl0sWzAsMiwiIiwxLHsic3R5bGUiOnsibmFtZSI6ImNvcm5lciJ9fV1d
\[\begin{tikzcd}
	K & A \\
	{C(P)} & {CS(A)}
	\arrow["{C(g)}"', from=2-1, to=2-2]
	\arrow["{\eta_A}", from=1-2, to=2-2]
	\arrow["{\pi_A}", from=1-1, to=1-2]
	\arrow["{\pi_{C(P)}}"', from=1-1, to=2-1]
	\arrow["\lrcorner"{anchor=center, pos=0.125}, draw=none, from=1-1, to=2-2]
\end{tikzcd}\]
is a pullback ($\eta_A$ is the component in $A$ of the unit of the adjunction $S \dashv C$). It can be proved that $S_A$ is the left adjoint of $C_A$.

\begin{definition}[\cite{janelidze1990pure}]
A relatively admissible adjunction $(S,C,\mathcal{A},\mathcal{P})$ is \emph{admissible} when the functor $C_A$ is full and faithful for every object $A$ of $\mathcal{A}$.
    
\end{definition}

Now, we recall the definition of effective descent morphism, as this concept will be of crucial importance for the detailed analysis of admissible Galois structures.

\begin{definition}

Let $\mathcal{A}$ be an admissible class of arrows in a category $\mathbb{A}$ with pullbacks. An arrow $h \colon A \rightarrow B$ is an \emph{effective descent} morphism relatively to $\mathcal{A}$ if:
\begin{itemize}
    \item $h \in \mathcal{A}$;
    \item the change-of-base functor $h^* \colon \mathcal{A}/B \rightarrow \mathcal{A}/A$ is monadic.
\end{itemize}

\end{definition}

It has been shown that, in Barr-exact categories, effective descent morphisms relatively to the class of all morphisms (that we will call, for simplicity, effective descent morphisms) are precisely the regular epimorphisms.\\

Finally, we observe how the notion of effective descent morphism enables the study of specific classes of arrows with respect to an admissible Galois structure.

\begin{definition}[\cite{janelidze1990pure}]
Given a relatively admissible adjunction $(S,C,\mathcal{A},\mathcal{P})$
\begin{itemize}
    \item a \emph{trivial extension} is an arrow $f \colon A \rightarrow B$ of $\mathcal{A}$ such that the square
% https://q.uiver.app/?q=WzAsNCxbMCwwLCJBIl0sWzAsMSwiQiJdLFsxLDAsIlBTKEEpIl0sWzEsMSwiUFMoQikiXSxbMCwxLCJmIiwyXSxbMiwzLCJQUyhmKSJdLFswLDIsIlxcZXRhX0EiXSxbMSwzLCJcXGV0YV9CIiwyXSxbMCwzLCIiLDEseyJzdHlsZSI6eyJuYW1lIjoiY29ybmVyIn19XV0=
\[\begin{tikzcd}
	A & {PS(A)} \\
	B & {PS(B)}
	\arrow["f"', from=1-1, to=2-1]
	\arrow["{PS(f)}", from=1-2, to=2-2]
	\arrow["{\eta_A}", from=1-1, to=1-2]
	\arrow["{\eta_B}"', from=2-1, to=2-2]
	\arrow["\lrcorner"{anchor=center, pos=0.125}, draw=none, from=1-1, to=2-2]
\end{tikzcd}\]
is a pullback (where $\eta$ is the unit of the adjunction $S \dashv P$);
\item a \emph{normal extension} is an arrow $f$ of $\mathcal{A}$ such that it is an effective descent morphism relatively to $\mathcal{A}$ and its kernel pair projections are trivial extensions;
\item a \emph{central extension} is an arrow $f$ of $\mathcal{A}$ such that there exists an effective descent morphism $g$ relatively to $\mathcal{A}$ and the pullback $g^*(f)$ of $f$ along $g$ is a trivial extension.
\end{itemize}

\end{definition}

Finally, we are ready to recall the notion of factorization system.

\begin{definition}
A \emph{factorization system} for a category $\cat{C}$ is a pair of classes of arrows $(\facte, \factm)$ such that:
\begin{itemize}
    \item for every $e \in \facte$ and $m \in \factm$ one has $e \downarrow m$, i.e.\ for every commutative square in $\cat{C}$
    % https://q.uiver.app/?q=WzAsNCxbMCwwLCJBIl0sWzEsMCwiQiJdLFswLDEsIkMiXSxbMSwxLCJEIl0sWzAsMSwiZSBcXGluIFxcZmFjdGUiXSxbMiwzLCJtIFxcaW4gXFxmYWN0bSJdLFswLDIsImciLDJdLFsxLDMsImgiXSxbMSwyLCJcXGV4aXN0cyFkIiwxLHsic3R5bGUiOnsiYm9keSI6eyJuYW1lIjoiZGFzaGVkIn19fV1d
\[\begin{tikzcd}
	A & B \\
	C & D
	\arrow["{e \in \facte}", from=1-1, to=1-2]
	\arrow["{m \in \factm}"', from=2-1, to=2-2]
	\arrow["g"', from=1-1, to=2-1]
	\arrow["h", from=1-2, to=2-2]
	\arrow["{\exists!d}"{description}, from=1-2, to=2-1]
\end{tikzcd}\]
there exists a unique arrow $d \colon B \rightarrow C$ such that $de=g$ and $md=h$; 
    \item every arrow $f$ in $\cat{C}$ factors as $f=me$, where $m \in \factm$ and $e \in \facte.$
\end{itemize}
A factorization system $(\facte, \factm)$ is \emph{stable} if the pullback of every arrow of $\facte$ along an arbitrary arrow is an arrow of $\facte$, too.
\end{definition}

\section{A Galois Theory for MV-Algebras}\label{Galois Theory for MV-Algebras}

In this section, we will study the category of MV-algebras from the perspective of categorical Galois theory. Given an MV-algebra A, we denote its radical, i.e. the intersection of its maximal ideals, by $\Rad(A)$. An MV-algebra is said to be \emph{semisimple} if its radical is trivial. It has been shown (for a proof of this fact see \cite{cignoli2013algebraic}) that an MV-algebra is semisimple if and only if it is a subdirect product of subalgebras of the MV-algebra [0,1]. We will also consider perfect MV-algebras, which are defined as follows: an MV-algebra A is said to be \emph{perfect} if $A = \Rad(A) \cup \lnot \Rad(A)$. We will show that the full subcategory of semisimple MV-algebras and the full subcategory of perfect MV-algebras constitute the torsion-free and torsion parts, respectively, of a pretorsion theory on $\mv$. We will then study the Galois theory defined by the reflector of the subcategory of semisimple MV-algebras.

\begin{definition}[\cite{cignoli2013algebraic}]

Given an MV-Algebra $A$, the \emph{radical} of A is defined as $$\Rad(A) \coloneqq \bigcap \{M \subseteq A \, | \, M \emph{ is a maximal ideal} \}.$$

\end{definition}

Let us consider an MV-algebra $A$. In \cite{cignoli2013algebraic}, Proposition 3.6.4 the authors show that $$\Rad(A)=\Inf(A) \cup \{ 0 \},$$ where $a \in \Inf(A)$ if and only if $a \neq 0$ and $na \leq \lnot a$ for every $n \in \mathbb{N}$. Moreover, in \cite{cignoli2013algebraic}, Lemma 7.3.3 it is proved that $$\Rad(A)=\bigvee \{J \subseteq A \, | \, J \text{ is a nilpotent ideal of } A \},$$ where an ideal $J$ is said to be nilpotent if, for every $x,y \in J$, one has $x \odot y=0$, and the join is computed in the poset $\Ideals(A)$ of ideals of $A$.

\begin{definition}[\cite{cignoli2013algebraic}]
Let $A$ be an MV-algebra. $A$ is \emph{semisimple} if its radical $\Rad(A)$ is trivial, i.e. $$\Rad(A) = \{0 \}.$$ 
$A$ is \emph{perfect} if it can be expressed as the union of its radical $\Rad(A)$ and the negation of its radical, i.e. $$A=\Rad(A) \cup \lnot \Rad(A),$$ where, for a subset $S$ of $A$, $\lnot S \coloneqq \{x \in A \, | \, \lnot x \in S \}$
\end{definition}

Given an MV-algebra $A$, we define $S(A) \coloneqq A/\Rad(A)$. Thanks to \cite{cignoli2013algebraic}, Lemma 3.6.6 we obtain that $S(A)$ is semisimple.

\begin{remark}

Given an MV-algebra $A$, we define $P(A) \coloneqq \Rad(A) \cup \lnot \Rad(A)$. $P(A)$ is a subalgebra of $A$:
\begin{itemize}
    \item $0 \in P(A)$;
    \item $x \in P(A)$ implies $\lnot x \in P(A)$;
    \item if $x,y \in P(A)$ then $x \oplus y \in P(A)$. To show this, we work on cases: if $x,y \in \Rad(A)$ then $x \oplus y \in \Rad(A)$; if $x \in \Rad(A)$ and $y \in \lnot \Rad(A)$ then, since $\lnot \Rad(A)$ is a filter and $y \leq x \oplus y$, we obtain $x \oplus y \in P(A)$; finally, if $x,y \in \lnot \Rad(A)$, since $\lnot \Rad(A)$ is a filter and $y \leq x \oplus y$, we get $x \oplus y \in P(A)$.
\end{itemize}
Moreover, it is easy to see that $P(A)$ is perfect.
\end{remark}
It is straightforward to verify that the argument just presented holds for any arbitrary ideal $I$ of $A$. In other words, the set $I \cup \lnot I$ always forms a subalgebra of $A$.\\
If we denote with $s\mv$ the full subcategory of $\mv$ whose objects are the semisimple MV-algebras, we get a functor $S$ described by the following assignment:
% https://q.uiver.app/?q=WzAsNixbMCwwLCJcXG12Il0sWzEsMCwic1xcbXYiXSxbMCwxLCJBIl0sWzEsMSwiUyhBKT1BL1xcUmFkKEEpIl0sWzAsNCwiQiJdLFsxLDQsIlMoQik9Qi9cXFJhZChBKSJdLFsyLDQsImYiLDJdLFsyLDMsIiIsMCx7InN0eWxlIjp7InRhaWwiOnsibmFtZSI6Im1hcHMgdG8ifX19XSxbNCw1LCIiLDIseyJzdHlsZSI6eyJ0YWlsIjp7Im5hbWUiOiJtYXBzIHRvIn19fV0sWzMsNSwiUyhmKT1cXGJhcntmfSJdLFswLDFdXQ==
\[\begin{tikzcd}[row sep=5pt]
	\mv & s\mv \\
	A & {S(A)=A/\Rad(A)} \\
	\\
	\\
	B & {S(B)=B/\Rad(B),}
	\arrow["S",from=1-1, to=1-2]
	\arrow["f"', from=2-1, to=5-1]
	\arrow[maps to, from=2-1, to=2-2]
	\arrow[maps to, from=5-1, to=5-2]
	\arrow["{S(f)=\overline{f}}", from=2-2, to=5-2]
	\arrow[from=1-1, to=1-2]
\end{tikzcd}\]
where $\overline{f}([a]) \coloneqq [f(a)]$, for every $[a] \in S(A)$. We have to show that $\overline{f}$ is well defined: $a \in \Rad(A)$ if and only if $na \leq \lnot a$ for every $n \in \mathbb{N}$, thus $nf(a)=f(na) \leq f(\lnot a)=\lnot f(a)$ for every $n \in \mathbb{N}$, and so $f(a) \in \Rad(B)$.\\
Similarly, if we denote with $p\mv$ the full subcategory of $\mv$ whose objects are the perfect MV-algebras, we get a functor $P$ described by the following assignment
% https://q.uiver.app/?q=WzAsNixbMCwwLCJcXG12Il0sWzEsMCwicFxcbXYiXSxbMCwxLCJBIl0sWzEsMSwiUChBKT1cXFJhZChBKSBcXGN1cCBcXGxub3QgXFxSYWQoQSkiXSxbMCw0LCJCIl0sWzEsNCwiUChCKT1cXFJhZChCKSBcXGN1cCBcXGxub3QgXFxSYWQoQSkiXSxbMiw0LCJmIiwyXSxbMiwzLCIiLDAseyJzdHlsZSI6eyJ0YWlsIjp7Im5hbWUiOiJtYXBzIHRvIn19fV0sWzQsNSwiIiwyLHsic3R5bGUiOnsidGFpbCI6eyJuYW1lIjoibWFwcyB0byJ9fX1dLFszLDUsIlAoZikiXSxbMCwxLCJQIl1d
\[\begin{tikzcd}[row sep=5pt]
	\mv & p\mv \\
	A & {P(A)=\Rad(A) \cup \lnot \Rad(A)} \\
	\\
	\\
	B & {P(B)=\Rad(B) \cup \lnot \Rad(A),}
	\arrow["f"', from=2-1, to=5-1]
	\arrow[maps to, from=2-1, to=2-2]
	\arrow[maps to, from=5-1, to=5-2]
	\arrow["{P(f)}", from=2-2, to=5-2]
	\arrow["P", from=1-1, to=1-2]
\end{tikzcd}\]
where $P(f)(a)=f(a)$ for every $a \in P(A).$ We have to show that $P(f)$ is well defined: as we saw before, if $a \in \Rad(A)$ then $f(a) \in \Rad(B)$, and so $f \colon A \rightarrow B$ restricts to $P(f) \colon P(A) \rightarrow P(B)$.

\begin{proposition}

The inclusion functor $i \colon s\mv \hookrightarrow \mv$ is right adjoint to $S$.

\begin{proof}

We construct the unit of the adjunction $\eta \colon id_{\mv} \rightarrow iS$ as the quotient projection $\eta_A \colon A \rightarrow A / \Rad(A)$, for every MV-algebra $A$. Clearly $\eta$ is a natural transformation. Moreover, $\eta$ satisfies the universal property of the unit: given a morphism of MV-algebras $g \colon A \rightarrow i(B)$, where $B$ is semisimple, we want to prove that there exists a unique morphism of MV-algebras $\tilde{g} \colon S(A) \rightarrow B$ such that the following diagram commutes:
% https://q.uiver.app/?q=WzAsMyxbMCwwLCJBIl0sWzIsMCwiaVMoQSkiXSxbMSwxLCJpKEIpIl0sWzAsMSwiXFxldGFfQSJdLFswLDIsImciLDJdLFsxLDIsImkoXFx0aWxkZXtnfSkiXV0=
\[\begin{tikzcd}
	A && {iS(A)} \\
	& {i(B).}
	\arrow["{\eta_A}", from=1-1, to=1-3]
	\arrow["g"', from=1-1, to=2-2]
	\arrow["{i(\tilde{g})}", from=1-3, to=2-2]
\end{tikzcd}\]
Since $B$ is semisimple, if $a \in \Rad(A)$ we have $g(a)=0$; therefore $g$ defines a unique morphism $\tilde{g}$ such that $\tilde{g}([a])=g(a)$, for every $a \in A$.
\end{proof}

\end{proposition}

\begin{proposition}

The inclusion functor $j \colon p\mv \hookrightarrow \mv$ is left adjoint to $P$.

\begin{proof}

We construct the counit of the adjunction $\varepsilon \colon jP \rightarrow id_{\mv}$ as the inclusion $\varepsilon_A \colon j(P(A)) \rightarrow A$, for every MV-algebra $A$. Clearly $\varepsilon$ is a natural transformation. Moreover, $\varepsilon$ satisfies the universal property of the counit: given a morphism of MV-algebras $h \colon j(B) \rightarrow A$, where $B$ is perfect, we want to prove that there exists a unique morphism of MV-algebras $\tilde{h} \colon B \rightarrow P(A)$ such that the following diagram commutes 
% https://q.uiver.app/?q=WzAsMyxbMSwwLCJqKEIpIl0sWzAsMSwiaihQKEEpKSJdLFsyLDEsIkEiXSxbMCwyLCJoIl0sWzEsMiwiXFx2YXJlcHNpbG9uX0EiLDJdLFswLDEsImooXFx0aWxkZXtofSkiLDJdXQ==
\[\begin{tikzcd}
	& {j(B)} \\
	{j(P(A))} && A.
	\arrow["h", from=1-2, to=2-3]
	\arrow["{\varepsilon_A}"', from=2-1, to=2-3]
	\arrow["{j(\tilde{h})}"', from=1-2, to=2-1]
\end{tikzcd}\]
Since $B$ is perfect, then, for every $b \in B$, $b \in \Rad(B)$ or $b \in \lnot \Rad(B)$ and $h(b) \in \Rad(A)$ or $h(b) \in \lnot \Rad(A)$, i.e.\ $h(b) \in P(A)$; therefore $h$ defines a unique morphism $\tilde{h}$ such that $\tilde{h}(b)=h(b)$, for every $b \in B$.
\end{proof}

\end{proposition}

More specifically, as anticipated at the beginning of this section, the pair of subcategories just introduced allows us to introduce a pretorsion theory in $\mv$. To show this, we define the class of trivial objects as $$\zeros \coloneqq \{ \terminal, \initial \},$$ where $\terminal$ indicates the terminal object of $\mv$ and $\initial$ the initial object.

\begin{proposition}

$(p\mv,s\mv)$ is a pretorsion theory for $\mv$.

\begin{proof}

First of all we observe that $$p\mv \cap s\mv= \{ \terminal, \initial \}.$$ We prove that, given two MV-algebras $A,B$, with $A$ perfect and $B$ semisimple, $\mv(A,B)\subseteq \zeroarrows$. Let us suppose $B \neq \terminal$ (this implies $A \neq \terminal$); then a morphism $f \colon A \rightarrow B$ factors in the following way:
% https://q.uiver.app/?q=WzAsMyxbMCwwLCJBIl0sWzIsMCwiQiJdLFsxLDEsIlxcY2F0ezJ9Il0sWzAsMSwiZiJdLFswLDIsIlxcQ2hpX3tcXGxub3QgXFxSYWQoQSl9IiwyXSxbMiwxLCIiLDIseyJzdHlsZSI6eyJ0YWlsIjp7Im5hbWUiOiJob29rIiwic2lkZSI6InRvcCJ9fX1dXQ==
\[\begin{tikzcd}
	A && B \\
	& {\initial,}
	\arrow["f", from=1-1, to=1-3]
	\arrow["{\chi_{\lnot \Rad(A)}}"', from=1-1, to=2-2]
	\arrow[hook, from=2-2, to=1-3]
\end{tikzcd}\]
where
\[ \chi_{\lnot \Rad(A)}(a) \coloneqq \begin{cases} 
          0 & a \in \Rad(A) \\
          1 & a \in \lnot \Rad(A).
       \end{cases}
    \]
    In fact, if $a \in A$ then $a \in \Rad(A)$ or $a \in \lnot \Rad(A)$, since $A$ is perfect, and so $f(a) \in \Rad(B)=\{ 0 \}$ or $f(a) \in \lnot \Rad(A)=\{ 1 \}$, because $B$ is semisimple.
If $B= \terminal$ we have 
% https://q.uiver.app/?q=WzAsMyxbMCwwLCJBIl0sWzIsMCwiXFxjYXR7MX0iXSxbMSwxLCJcXGNhdHsxfSJdLFswLDEsImYiXSxbMCwyLCIiLDIseyJzdHlsZSI6eyJib2R5Ijp7Im5hbWUiOiJzcXVpZ2dseSJ9fX1dXQ==
\[\begin{tikzcd}
	A && {\terminal} \\
	& {\terminal.}
	\arrow["f", from=1-1, to=1-3]
	\arrow["f"', from=1-1, to=2-2]
	\arrow[equals, from=2-2, to=1-3]
\end{tikzcd}\]
Now, we want to show that, for every MV-algebra $B$, we have a pre-exact sequence with $B$ in the middle, a torsion object on the left, and a torsion-free object on the right. If $B= \terminal$ the sequence is given by $B=B=B$. Suppose $B \neq \terminal$; we prove that 
% https://q.uiver.app/?q=WzAsMyxbMCwwLCJQKEIpIl0sWzEsMCwiQiJdLFsyLDAsIlMoQikiXSxbMCwxLCJcXHZhcmVwc2lsb25fQiIsMix7InN0eWxlIjp7InRhaWwiOnsibmFtZSI6Imhvb2siLCJzaWRlIjoidG9wIn19fV0sWzEsMiwiXFxldGFfQiIsMix7InN0eWxlIjp7ImhlYWQiOnsibmFtZSI6ImVwaSJ9fX1dXQ==
\[\begin{tikzcd}
	{P(B)} & B & {S(B)}
	\arrow["{\varepsilon_B}"', hook, from=1-1, to=1-2]
	\arrow["{\eta_B}"', two heads, from=1-2, to=1-3]
\end{tikzcd}\]
is a pre-exact sequence. Observe that, since $\eta_B(b)=0$ for every $b \in \Rad(B)$, $\eta_B \varepsilon_B \in \zeroarrows$. We show that $\varepsilon_B$ is a prekernel of $\eta_B$. Consider an arrow $\lambda_B \colon C \rightarrow B$ such that there exists an arrow $x \colon C \rightarrow \initial$ making the diagram below commutative 
% https://q.uiver.app/?q=WzAsNCxbMCwwLCJCIl0sWzEsMCwiUyhCKSJdLFswLDEsIkMiXSxbMSwxLCJcXGNhdHsyfSJdLFswLDEsIlxcZXRhX0IiLDAseyJzdHlsZSI6eyJoZWFkIjp7Im5hbWUiOiJlcGkifX19XSxbMiwwLCJcXGxhbWJkYV9CIl0sWzMsMSwiIiwyLHsic3R5bGUiOnsidGFpbCI6eyJuYW1lIjoiaG9vayIsInNpZGUiOiJ0b3AifX19XSxbMiwzLCJ4IiwyXV0=
\[\begin{tikzcd}
	B & {S(B)} \\
	C & {\initial.}
	\arrow["{\eta_B}", two heads, from=1-1, to=1-2]
	\arrow["{\lambda_B}", from=2-1, to=1-1]
	\arrow[hook, from=2-2, to=1-2]
	\arrow["x"', from=2-1, to=2-2]
\end{tikzcd}\]
Notice that if $\eta_B \lambda_B$ factored through $\terminal$, then we would have $S(B)=\terminal$; since $S(B)=\terminal$ implies $1 \in \Rad(B)$, it would follow $B=\terminal$. Now, for every $c \in C$, if $x(c)=0$ then $\lambda_B(c) \in \Rad(B)$ and, if $x(c)=1$, then $\lambda_B(c) \in \lnot \Rad(B)$; therefore, $\lambda_B$ restricts to $\lambda'_B \colon C \rightarrow P(B)$ and the following diagram commutes
% https://q.uiver.app/?q=WzAsNSxbMCwwLCJQKEIpIl0sWzIsMCwiQiJdLFs0LDAsIlMoQikiXSxbMSwxLCJDIl0sWzMsMSwiXFxjYXR7Mn0iXSxbNCwyLCIiLDIseyJzdHlsZSI6eyJ0YWlsIjp7Im5hbWUiOiJob29rIiwic2lkZSI6InRvcCJ9fX1dLFszLDQsIngiLDJdLFszLDAsIlxcbGFtYmRhJ19CIl0sWzMsMSwiXFxsYW1iZGFfQiJdLFsxLDIsIlxcZXRhX0IiLDAseyJzdHlsZSI6eyJoZWFkIjp7Im5hbWUiOiJlcGkifX19XSxbMCwxLCJcXHZhcmVwc2lsb25fQiIsMCx7InN0eWxlIjp7InRhaWwiOnsibmFtZSI6Imhvb2siLCJzaWRlIjoidG9wIn19fV1d
\[\begin{tikzcd}[column sep=5pt]
	{P(B)} && B && {S(B)} \\
	& C && {\initial.}
	\arrow[hook, from=2-4, to=1-5]
	\arrow["x"', from=2-2, to=2-4]
	\arrow["{\lambda'_B}", from=2-2, to=1-1]
	\arrow["{\lambda_B}", from=2-2, to=1-3]
	\arrow["{\eta_B}", two heads, from=1-3, to=1-5]
	\arrow["{\varepsilon_B}", hook, from=1-1, to=1-3]
\end{tikzcd}\]
Therefore, $\varepsilon_B$ is a prekernel of $\eta_B$.\\
Next, we focus on $\eta_B$: we want to show that it is a precokernel of $\varepsilon_B$. Fix an arrow $\theta_B \colon B \rightarrow C$ such that $\theta_B \varepsilon_B \in \zeroarrows$. If $\theta_B \varepsilon_B$ factors through $\terminal$, we can conclude that $C=\terminal$. Hence, the claim becomes trivial. Then, suppose there exists an arrow $y \colon P(B) \rightarrow \initial$ making the following diagram commutes:
% https://q.uiver.app/?q=WzAsNCxbMCwwLCJQKEIpIl0sWzEsMCwiQiJdLFswLDEsIlxcY2F0ezJ9Il0sWzEsMSwiQyJdLFsxLDMsIlxcdGhldGFfQiJdLFswLDEsIlxcdmFyZXBzaWxvbl9CIiwwLHsic3R5bGUiOnsidGFpbCI6eyJuYW1lIjoiaG9vayIsInNpZGUiOiJ0b3AifX19XSxbMCwyLCJ5IiwyXSxbMiwzLCIiLDIseyJzdHlsZSI6eyJ0YWlsIjp7Im5hbWUiOiJob29rIiwic2lkZSI6InRvcCJ9fX1dXQ==
\[\begin{tikzcd}
	{P(B)} & B \\
	{\initial} & C.
	\arrow["{\theta_B}", from=1-2, to=2-2]
	\arrow["{\varepsilon_B}", hook, from=1-1, to=1-2]
	\arrow["y"', from=1-1, to=2-1]
	\arrow[hook, from=2-1, to=2-2]
\end{tikzcd}\]
If $b \in \Rad(B)$ then $\theta_B(b) \in \initial \subseteq C$; but, if $\theta_B(b)=1$, we get $1 \in \Rad(C)$ and so $C=\terminal$ (we are excluding this case). Hence, $\theta_B$ defines a unique morphism $\theta'_B \colon S(B) \rightarrow C$ such that the diagram below is commutative
% https://q.uiver.app/?q=WzAsNSxbMCwwLCJQKEIpIl0sWzIsMCwiQiJdLFs0LDAsIlMoQikiXSxbMywxLCJDIl0sWzEsMSwiXFxjYXR7Mn0iXSxbMCwxLCJcXHZhcmVwc2lsb25fQiIsMCx7InN0eWxlIjp7InRhaWwiOnsibmFtZSI6Imhvb2siLCJzaWRlIjoidG9wIn19fV0sWzEsMiwiXFxldGFfQiIsMCx7InN0eWxlIjp7ImhlYWQiOnsibmFtZSI6ImVwaSJ9fX1dLFs0LDMsIiIsMCx7InN0eWxlIjp7InRhaWwiOnsibmFtZSI6Imhvb2siLCJzaWRlIjoidG9wIn19fV0sWzAsNCwieSIsMl0sWzEsMywiXFx0aGV0YV9CIl0sWzIsMywiXFx0aGV0YSdfQiJdXQ==
\[\begin{tikzcd}[column sep=5pt]
	{P(B)} && B && {S(B)} \\
	& {\initial} && C.
	\arrow["{\varepsilon_B}", hook, from=1-1, to=1-3]
	\arrow["{\eta_B}", two heads, from=1-3, to=1-5]
	\arrow[hook, from=2-2, to=2-4]
	\arrow["y"', from=1-1, to=2-2]
	\arrow["{\theta_B}", from=1-3, to=2-4]
	\arrow["{\theta'_B}", from=1-5, to=2-4]
\end{tikzcd}\]
Therefore, $\eta_B$ is a precokernel of $\varepsilon_B$.
\end{proof}

\end{proposition}
The purpose of the following part of the section is to examine the adjunction $$S \dashv i \colon s \mv \hookrightarrow \mv$$ from the perspective of categorical Galois theory.\\
In the literature, adjunctions studied through the lens of categorical Galois theory are frequently adjunctions between semi-abelian categories. Additionally, in cases where the adjunction is induced by the reflection of a reflective subcategory, it is often required that the subcategory is a Birkhoff subcategory. A full and reflective subcategory $\cat{D}$ of a category $\cat{C}$ is considered a \emph{Birkhoff} subcategory if it is closed under subobjects and regular images.
It is important to mention that our work does not fall under the conditions specified above. In fact, our category $\mv$ is not semi-abelian (it is not pointed), and the full subcategory $s \mv \hookrightarrow \mv$ is not a Birkhoff subcategory. In fact, let us consider the MV-algebra $$A \coloneqq \prod_{n \geq 2} L_n$$ where $L_n \coloneqq \{0,\frac{1}{n}, \dots, \frac{n-1}{n},1 \} \subseteq [0,1]$ and the operations are the ones induced by $[0,1]$. First of all, let us show that $A$ is semisimple. Given an element $x=(x_n)_{n\geq2} \in A$, with $x \neq 0$, there exists a natural number $n \geq 2$ and a natural number $0 < k \leq n$ such that $x_n = \frac{k}{n} \neq 0$, hence $(\lnot x)_n=\lnot x_n=1-\frac{k}{n}=\frac{n-k}{n} \neq 1$. Therefore, there must exist a natural number $m$ such that $mk>n-k$ and so $mx \nleq \lnot x$, i.e.\ $\Inf(A)= \emptyset$ and $\Rad(A)=\{ 0 \}$. Next, we consider the relation $\rho \subseteq A \times A$ defined by $(x,y) \in \rho$ if and only if there exists a natural number $N \geq 2$ such that $x_n=y_n$ for every $n \geq N$, where $x=(x_n)_{n \geq 2}$ and $y=(y_n)_{n \geq 2}$. It is clear that $\rho$ is a congruence. Let $[z] \in A/\rho$ where $z=(z_n)_{n \geq 2}$ and $z_n=\frac{1}{n}$ for every $n \geq 2$. Our goal now is to show that $[z] \in \Inf(A/\rho)$. Fix a natural number $m \in \mathbb{N}$; clearly $\frac{m}{n} \leq \frac{n-1}{n}$ for every $n>m$, that is $mz_n \leq \lnot z_n$. We define $$y_n \coloneqq 
\begin{cases}
0 \text{ if } n \leq m \\
z_n \text{ if } n>m.
\end{cases}$$
Then $m[z]=m[y] \leq \lnot [y]=\lnot [z]$, and so the radical $\Rad(A/\rho)$ is not trivial (i.e.\ $A/\rho$ is not semisimple).

\begin{proposition}

For every MV-algebra $B$ the counit of the adjunction $$S_B \dashv i_B \colon s\mv/s(B) \rightarrow \mv/B$$ is an isomorphism.

\begin{proof}

When $B = \terminal$ the assertion is trivial. Suppose $B \neq \terminal$; we recall that, for every object $f \colon A \rightarrow B$ of $\mv/B$, $S_B(f \colon A \rightarrow B)=(\overline{f} \colon S(A) \rightarrow S(B))$ and, for every object $\varphi \colon A \rightarrow S(B)$ of $s\mv/s(B)$, $i_B(\varphi \colon A \rightarrow S(B))=(\varphi' \colon A' \rightarrow B)$ is defined by the following pullback:
% https://q.uiver.app/?q=WzAsNCxbMCwwLCJBJyJdLFsxLDAsIkEiXSxbMSwxLCJTKEIpPUIvXFxyYWQoQikiXSxbMCwxLCJCIl0sWzMsMiwiXFxldGFfQiIsMix7InN0eWxlIjp7ImhlYWQiOnsibmFtZSI6ImVwaSJ9fX1dLFsxLDIsIlxcdmFycGhpIiwwLHsic3R5bGUiOnsiaGVhZCI6eyJuYW1lIjoiZXBpIn19fV0sWzAsMywiXFx2YXJwaGknIiwyLHsic3R5bGUiOnsiaGVhZCI6eyJuYW1lIjoiZXBpIn19fV0sWzAsMSwiXFx2YXJwaGknJyIsMCx7InN0eWxlIjp7ImhlYWQiOnsibmFtZSI6ImVwaSJ9fX1dLFswLDIsIiIsMSx7InN0eWxlIjp7Im5hbWUiOiJjb3JuZXIifX1dXQ==
\[\begin{tikzcd}
	{A'} & A \\
	B & {S(B).}
	\arrow["{\eta_B}"', two heads, from=2-1, to=2-2]
	\arrow["\varphi", from=1-2, to=2-2]
	\arrow["{\varphi'}"', from=1-1, to=2-1]
	\arrow["{\varphi''}", two heads, from=1-1, to=1-2]
	\arrow["\lrcorner"{anchor=center, pos=0.125}, draw=none, from=1-1, to=2-2]
\end{tikzcd}\]
We want to prove that $$(S(\varphi') \colon S(A') \rightarrow S(B)) \cong (\varphi \colon A \rightarrow S(B)).$$
Consider the following commutative diagram, where $\varphi'''$ is induced by the universal property of $\eta$,
% https://q.uiver.app/?q=WzAsMyxbMCwwLCJBJyJdLFsyLDAsIkEnL1xcUmFkKEEnKSJdLFsxLDEsIkEiXSxbMCwxLCJcXGV0YV97QSd9IiwwLHsic3R5bGUiOnsiaGVhZCI6eyJuYW1lIjoiZXBpIn19fV0sWzAsMiwiXFx2YXJwaGknJyIsMix7InN0eWxlIjp7ImhlYWQiOnsibmFtZSI6ImVwaSJ9fX1dLFsxLDIsIlxcdmFycGhpJycnIiwwLHsic3R5bGUiOnsiaGVhZCI6eyJuYW1lIjoiZXBpIn19fV1d
\[\begin{tikzcd}[column sep=5pt]
	{A'} && {A'/\Rad(A')} \\
	& A.
	\arrow["{\eta_{A'}}", two heads, from=1-1, to=1-3]
	\arrow["{\varphi''}"', two heads, from=1-1, to=2-2]
	\arrow["{\varphi'''}", two heads, from=1-3, to=2-2]
\end{tikzcd}\]
Clearly $\varphi'''$ is surjective; let us show that it is also injective. Fix an element $(b,a) \in A'$, where $a \in A$, $b \in B$ and $\varphi(a)=\eta_B(b)=[b]$. If $\varphi'''([b,a])=0$, we obtain $\varphi'''\eta_{A'}(b,a)=0$ and so $\varphi''(b,a)=0$, i.e.\ $a=0$. We want 
to prove that $(b,0) \in \Rad(A')$: since $b \in \Rad(B)$ we know that $nb \leq \lnot b$ for every $n \in \mathbb{N}$; then $n(b,0)=(nb,0) \leq (\lnot b, 1)= \lnot (b,0)$, for every $n \in \mathbb{N}$, and so $(b,0) \in \Rad(A')$. Therefore $\kerr(\varphi''')$ is trivial, i.e.\ $\varphi'''$ is injective. It remains to prove that the following diagram is commutative:
% https://q.uiver.app/?q=WzAsMyxbMCwwLCJBJy9cXFJhZChBJykiXSxbMSwxLCJCL1xcUmFkKEIpIl0sWzAsMiwiQSJdLFswLDIsIlxcdmFycGhpJycnIiwyXSxbMCwxLCJTKFxcdmFycGhpJykiXSxbMiwxLCJcXHZhcnBoaSIsMl1d
\[\begin{tikzcd}[row sep=5pt]
	{A'/\Rad(A')} \\
	& {B/\Rad(B)} \\
	A.
	\arrow["{\varphi'''}"', from=1-1, to=3-1]
	\arrow["{S(\varphi')}", from=1-1, to=2-2]
	\arrow["\varphi"', from=3-1, to=2-2]
\end{tikzcd}\]
Take an element $(b,a) \in A'$. We know that $\varphi \varphi'''([b,a])=\varphi(a)$, $\varphi(a)=[b]$, and $S(\varphi')([b,a])=[b]$.
\end{proof}

\end{proposition}

\begin{comment}

\begin{proposition}

Finite limits in $s\mv$ are computed as in $\mv$.

\begin{proof}

We start dealing with binary products. Let us consider two objects $A,B$ of $s\mv$. We consider an element $(a,b) \in \Rad(A \times B)$, therefore for each $n \in \mathbb{N}$ we have $n(a,b) \leq \lnot (a,b)$ and so $na \leq \lnot a$, $nb \leq \lnot b$ i.e.\ $a \in \Rad(A)=\{ 0 \}$ and $b \in \Rad(B)= \{ 0 \}$. Then, we obtain $\Rad(A \times B)= \{ 0 \}$. Now, let us study the equalizers in $s \mv$. We fix two arrows $f,g \colon A \rightarrow B$ of $s\mv$. We know that the equalizer of $f,g$ in $\mv$ is defined by $E \hookrightarrow A$ where $E \coloneqq \{a \in A\, | \, f(a)=g(g) \}$. If we show that $E$ is semisimple  we are done. But $\Rad(E)=\{ a \in A \, | \, na \leq \lnot \text{ for every } n \in \mathbb{N}, a \in E \}=\Rad(A) \cap E= \{ 0 \}$ since $A$ is semisimple.

\end{proof}

\end{proposition}

\end{comment}

We observe that, thanks to the previous proposition, the adjunction $S \dashv i$ is admissible relatively to the class of all arrows in $\mv$ and the class of all arrows in $s\mv$. We will denote this Galois structure by $\Gamma$.

\begin{proposition}\label{proposizione da sostituire}

Given an arrow $f \colon A \rightarrow B$ in $\mv$, the restriction of $f$ (considered as map) $f \colon \Rad(A) \rightarrow \Rad(B)$ is injective if and only if $P(f) \colon P(A) \rightarrow P(B)$ is injective, or $P(A)= \initial$ and $P(B)= \terminal$.

\begin{proof}

Thanks to Lemma \ref{ideali}, we know that $f \colon \Rad(A) \rightarrow \Rad(B)$ is injective if and only if $\kerr(f) \cap \Rad(A) = \{ 0 \}$. Suppose first that $B \neq \terminal$ (which implies $A \neq \terminal$). If $a \in \kerr(P(f))$, then we have $a \in \Rad(A)$: indeed, if $a \in \lnot \Rad(A)$, then $f(a) = 0 \in \lnot \Rad(B)$, which contradicts our assumption that $B \neq \terminal$. Hence, we have $\kerr(P(f))=\kerr(f) \cap \Rad(A)$, and so $P(f)$ is injective if and only if $\kerr(f) \cap \Rad(A)=\{ 0 \}$. Now suppose that $B = \terminal$. Then  the restriction $f \colon \Rad(A) \rightarrow \Rad(B)=\{ 0 \}$ is injective if and only if $\Rad(A)= \{ 0 \}$. Indeed, the condition $\Rad(A)= \{ 0 \}$ is equivalent to either $P(A)=\initial$ or $P(A)=\terminal$.
\end{proof}

\end{proposition}

\begin{proposition}

Given an arrow $f \colon A \rightarrow B$ in $\mv$, the following diagram is a pullback (i.e.\ $f$ is a trivial extension for $\Gamma$)
% https://q.uiver.app/?q=WzAsNCxbMCwwLCJBIl0sWzAsMSwiQiJdLFsxLDAsIlMoQSk9QS9cXFJhZChBKSJdLFsxLDEsIlMoQik9Qi9cXFJhZChCKSJdLFswLDIsIlxcZXRhX0EiLDAseyJzdHlsZSI6eyJoZWFkIjp7Im5hbWUiOiJlcGkifX19XSxbMSwzLCJcXGV0YV9CIiwyLHsic3R5bGUiOnsiaGVhZCI6eyJuYW1lIjoiZXBpIn19fV0sWzAsMSwiZiIsMix7InN0eWxlIjp7ImhlYWQiOnsibmFtZSI6ImVwaSJ9fX1dLFsyLDMsIlMoZik9XFxiYXJ7Zn0iXV0=
\[\begin{tikzcd}
	A & {S(A)=A/\Rad(A)} \\
	B & {S(B)=B/\Rad(B)}
	\arrow["{\eta_A}", two heads, from=1-1, to=1-2]
	\arrow["{\eta_B}"', two heads, from=2-1, to=2-2]
	\arrow["f"', from=1-1, to=2-1]
	\arrow["{S(f)=\overline{f}}", from=1-2, to=2-2]
\end{tikzcd}\]
if and only if $P(f) \colon P(A) \rightarrow P(B)$ is an isomorphism, or $P(A)= \initial$ and $P(B)=\terminal$.

\begin{proof}

We apply Lemma \ref{Lemma dei pullback} and we get that the above diagram is a pullback if and only if the restriction $f \colon \Rad(A) \rightarrow \Rad(B)$ is bijective. Thanks to Proposition \ref{proposizione da sostituire}, this last statement holds if and only if $P(f) \colon P(A) \rightarrow P(B)$ is an isomorphism, or $P(A)= \initial$ and $P(B)=\terminal$.
\end{proof}

\end{proposition}

\begin{proposition}

Given an effective descent morphism (i.e.\ a regular epimorphism) $f \colon A \twoheadrightarrow B$ in $\mv$, $f$ is a normal extension for $\Gamma$ if and only if $P(f) \colon P(A) \rightarrow P(B)$ is injective, or $P(A)= \initial$ and $P(B)= \terminal$.

\begin{proof}

Consider the kernel pair of $f$ defined by the pullback 
% https://q.uiver.app/?q=WzAsNCxbMCwxLCJBIl0sWzEsMSwiQiJdLFsxLDAsIkEiXSxbMCwwLCJLIl0sWzAsMSwiZiIsMix7InN0eWxlIjp7ImhlYWQiOnsibmFtZSI6ImVwaSJ9fX1dLFsyLDEsImYiLDAseyJzdHlsZSI6eyJoZWFkIjp7Im5hbWUiOiJlcGkifX19XSxbMywyLCJcXHBpXzIiLDAseyJzdHlsZSI6eyJoZWFkIjp7Im5hbWUiOiJlcGkifX19XSxbMywwLCJcXHBpXzEiLDIseyJzdHlsZSI6eyJoZWFkIjp7Im5hbWUiOiJlcGkifX19XSxbMywxLCIiLDEseyJzdHlsZSI6eyJuYW1lIjoiY29ybmVyIn19XV0=
\[\begin{tikzcd}
	\Eq(f) & A \\
	A & B.
	\arrow["f"', two heads, from=2-1, to=2-2]
	\arrow["f", two heads, from=1-2, to=2-2]
	\arrow["{\pi_2}", two heads, from=1-1, to=1-2]
	\arrow["{\pi_1}"', two heads, from=1-1, to=2-1]
	\arrow["\lrcorner"{anchor=center, pos=0.125}, draw=none, from=1-1, to=2-2]
\end{tikzcd}\]
We want to prove that $P(\pi_1)$ is an isomorphism if and only if $P(f)$ is injective, or $P(A) = \initial$ and $P(B)= \terminal$.
We start by assuming $B \neq \terminal$. Suppose that $P(f) \colon P(A) \rightarrow P(B)$ is injective; we observe that $P(\Eq(f))=\{ (a_1,a_2) \in A \times A \, | \, a_1,a_2 \in \Rad(A) \text{ and } f(a_1)=f(a_2) \} \cup \{ (a_1,a_2) \in A \times A \, | \, a_1,a_2 \in \lnot \Rad(A) \text{ and } f(a_1)=f(a_2) \}$. Therefore, since $P(f)$ is injective and $\Rad(A) \cap \lnot \Rad(A)= \emptyset$ (otherwise we would get $A= \terminal$), if $(a_1,a_2) \in P(\Eq(f))$ then $f(a_1)=f(a_2)$ and so $a_1=a_2$. Hence, $P(\Eq(f))=\{ (a,a) \in A \times A \, | \, a \in P(A) \}$ and clearly $P(\pi_1) \colon P(\Eq(f)) \rightarrow P(A)$ is an isomorphism. Conversely, if we assume that $P(\pi_1)$ is an isomorphism and we consider an element $a \in P(A)$ such that $P(f)(a)=0$, then $a \in \Rad(A)$ (otherwise we would obtain $B= \terminal$), and so $(0,a) \in P(\Eq(f))$ (in fact, $n(0,a)=(0,na) \leq (1,\lnot a)=\lnot (0,a)$ for every $n \in \mathbb{N}$). So, since $P(\pi_1)$ is an isomorphism and $P(\pi_1)(0,0)=P(\pi_1)(0,a)$, we get $a=0$ and then we deduce that $P(f)$ is injective. Finally, let us handle the case $B = \terminal$. If $A = \terminal$, then $\Eq(f)= \terminal$ and so the assertion is trivial. If $A \neq \terminal$ and $B = \terminal$ we have two possible situations. If $P(A)= \initial$, then $P(\Eq(f))=\initial$ and so $P(\pi_1)$ is an isomorphism and $P(f) \colon \initial \rightarrow \terminal$. If $P(A) \neq \initial$, hence $P(f) \colon P(A) \rightarrow \terminal$ is not injective and $P(\pi_1)$ is not an isomorphism. In fact, by assumption, we have an element $a \neq 0$ such that $a \in \Rad(A)$; therefore, we observe that there exist two different elements $(0,0), (0,a) \in P(\Eq(f))$ (since $f(a) \in \Rad(B)=\{ 0 \}$) such that $P(\pi_1)(0,0)=P(\pi_1)(0,a)$, and so $P(\pi_1)$ is not an isomorphism. 
\end{proof}

\end{proposition}

\begin{proposition}\label{proposizione centrali}

Given an effective descent morphism (i.e.\ a regular epimorphism) $f \colon A \twoheadrightarrow B$ in $\mv$, $f$ is a central extension for $\Gamma$ if and only if it is a normal extension, i.e.\ $P(f) \colon P(A) \rightarrow P(B)$ is injective, or $P(A)= \initial$ and $P(B)= \terminal$.

\begin{proof}

Consider the following pullback diagram:
% https://q.uiver.app/?q=WzAsNCxbMCwwLCJBIFxcdGltZXNfQiBDIl0sWzEsMCwiQyJdLFsxLDEsIkIiXSxbMCwxLCJBIl0sWzMsMiwiZiIsMix7InN0eWxlIjp7ImhlYWQiOnsibmFtZSI6ImVwaSJ9fX1dLFswLDMsIlxccGlfQSIsMix7InN0eWxlIjp7ImhlYWQiOnsibmFtZSI6ImVwaSJ9fX1dLFswLDEsIlxccGlfQyIsMCx7InN0eWxlIjp7ImhlYWQiOnsibmFtZSI6ImVwaSJ9fX1dLFsxLDIsImciLDAseyJzdHlsZSI6eyJoZWFkIjp7Im5hbWUiOiJlcGkifX19XSxbMCwyLCIiLDEseyJzdHlsZSI6eyJuYW1lIjoiY29ybmVyIn19XV0=
\[\begin{tikzcd}
	{A \times_B C} & C \\
	A & B.
	\arrow["f"', two heads, from=2-1, to=2-2]
	\arrow["{\pi_A}"', two heads, from=1-1, to=2-1]
	\arrow["{\pi_C}", two heads, from=1-1, to=1-2]
	\arrow["g", two heads, from=1-2, to=2-2]
	\arrow["\lrcorner"{anchor=center, pos=0.125}, draw=none, from=1-1, to=2-2]
\end{tikzcd}\]
Thanks to the previous proposition, if $P(f)$ is injective or $P(A)=\initial$ and $P(B)=\terminal$, then we can choose $g=f$. Conversely, suppose that there exists a regular epimorphism $g \colon C \twoheadrightarrow B$ such that the restriction of $\pi_C$ (considered as a map) $\pi_C \colon \Rad(A \times_B C) \rightarrow \Rad(C)$ is a bijection (i.e.\ $\pi_C$ is a trivial extension). We will show that the restriction of $f$ (considered as a map) $f \colon \Rad(A) \rightarrow \Rad(B)$ is injective. Thanks to Proposition \ref{proposizione da sostituire}, this is precisely equivalent to stating that $P(f)$ is injective, or $P(A)=\initial$ and $P(B)= \terminal$. So, it suffices to prove that $\kerr(f) \cap \Rad(A)=\{0\}$. To show this, let us fix an element $a \in \kerr(f) \cap \Rad(A)$. Then, $f(a)=0=g(0)$, and thus $(a,0)$ belongs to $A \times_B C$. Moreover, we observe that $(a,0)$ is an element of $\Rad(A \times_B C)$, since $a \in \Rad(A)$. Therefore, recalling that, by assumption, $\pi_C$ restricted to $\Rad(A \times_B C)$ is an injective map and observing that $\pi_C(a,0)=\pi_C(0,0)$, we deduce that $a=0$. Hence, the restriction of $f$ to $\Rad(A)$ is injective.
\end{proof}

\end{proposition}

We provide an example of a central extension which is not trivial. Consider the MV-algebra $$A \coloneqq \prod_{n \geq 2} L_n$$ and the congruence $\rho \subseteq A \times A$, introduced before, defined by $(x,y) \in \rho$ if and only if there exists a natural number $N \geq 2$ such that $x_n=y_n$ for every $n \geq N$, where $x=(x_n)_{n \geq 2}$ and $y=(y_n)_{n \geq 2}$. We prove that the quotient projection $$\pi \colon A \twoheadrightarrow A/\rho$$ is central but not trivial. Recalling that $A$ is semisimple, we immediately get that $P(A)= \initial$ and so the map $$P(\pi) \colon P(A)=\initial \rightarrow P(A/\rho)$$ is injective (since $A/\rho \neq \terminal$ implies $P(A/\rho) \neq \terminal$). However, the map $P(\pi)$ is not surjective. To see this, let $[z] \in A/\rho$ be fixed, where $z=(z_n)_{n \geq 2}$ and $z_n=\frac{1}{n}$ for every $n \geq 2$. We have already shown that $[z] \in \Inf(A/\rho)$, and it is clear that $[z]$ is not the zero element of $A/\rho$. Therefore, $[z]$ does not lie in the image of $P(\pi)$, and we conclude that $P(\pi)$ is not surjective.
This implies that the morphism $\pi$ is central but not trivial.\\

In \cite{everaert2011galois}, the authors show that, in the case of an adjunction between semi-abelian categories induced by the reflection of a Birkhoff subcategory, the Galois theory determined by such an adjunction allows for the introduction of a notion of \emph{commutator}. The goal of this part of the section is to present a construction similar to that of \cite{everaert2011galois} for the case of MV-algebras. Specifically, given a regular epimorphism $f \colon A \twoheadrightarrow B$, we seek to identify a subalgebra of $A$ such that $f$ is a central morphism if and only if the subalgebra is trivial (which in our case means that it is an element of the class $\zeros$).\\

Recalling that every MV-algebra is a distributive lattice with respect to the operations defined in Proposition \ref{inf}, the result presented in Proposition \ref{proposizione centrali} can be expressed in a different form: a regular epimorphism $f$ is a central extension if and only if $\kerr(f) \subseteq \pol{\Rad(A)}$. Here, for an MV-algebra $A$ and a non-empty subset $S \subseteq A$, we define the set $$\pol{S} \coloneqq \{ x \in A \, | \, x \land s=0 \text{ for every } s \in S \}.$$ It can be proved that $\pol{S}$ is an ideal of $A$.

\begin{lemma}

Consider a morphism $f \colon A \rightarrow B$ in $\mv$. Then, $\kerr(f) \cap \Rad(A)= \{ 0 \}$ if and only if $K[f] \cap P(A) \in \zeros$, where $K[f] \coloneqq A$ if $B = \terminal$, otherwise $K[f]$ is given by the following pullback:
% https://q.uiver.app/?q=WzAsNCxbMCwwLCJLW2ZdIl0sWzEsMCwiXFxpbml0aWFsIl0sWzEsMSwiQiJdLFswLDEsIkEiXSxbMywyLCJmIiwyXSxbMCwzLCJrIiwyXSxbMCwxXSxbMSwyLCJcXGliYW5ne0J9Il0sWzAsMiwiIiwxLHsic3R5bGUiOnsibmFtZSI6ImNvcm5lciJ9fV1d
\[\begin{tikzcd}
	{K[f]} & \initial \\
	A & B.
	\arrow["f"', from=2-1, to=2-2]
	\arrow["k"', from=1-1, to=2-1]
	\arrow[from=1-1, to=1-2]
	\arrow["{\ibang{B}}", from=1-2, to=2-2]
	\arrow["\lrcorner"{anchor=center, pos=0.125}, draw=none, from=1-1, to=2-2]
\end{tikzcd}\]

\begin{proof}
We notice that $K[f]$ is precisely given by the union of $I$ with $\lnot I$, where $I=\kerr(f)$. Therefore, we obtain $$\Rad(K[f])=\Rad(A) \cap K[f]= (\Rad(A) \cap \kerr(f)) \cup (\Rad(A) \cap \lnot \kerr(f)).$$
If $B \neq \terminal$, then $1 \notin \Rad(B)$ and so $\Rad(A) \cap \lnot \kerr(f)= \emptyset$. Hence $\kerr(f) \cap \Rad(A)= \{ 0 \}$ if and only if $\Rad(K[f])= \{ 0 \}$ (which means $P(K[f]) \in \zeros$). If $B = \terminal$, then $\kerr(f) = A$ and $K[f] = A$. Hence, the equivalence we need to prove simplifies to: $\Rad(A) = \{0 \}$ if and only if $P(A) \in \zeros$; this equivalence always holds, so the statement is trivial in this case.
\end{proof}

\end{lemma}
To conclude, based on what has been proven so far, we observe that a regular epimorphism $f \colon A \twoheadrightarrow B$ is central for $\Gamma$ if and only if $K[f] \cap P(A) \in \zeros$. It therefore makes sense to define, for a general regular epimorphism $f$, the following subalgebra of $A$: $$[A,K[f]]_{s\mv} \coloneqq K(f) \cap P(A);$$ $[A,K[f]]_{s\mv}$ has the following property: it belongs to $\zeros$ if and only if $f$ is central.\\

In the final part of this section, we will focus on the study of the functor $S \colon \mv \rightarrow s\mv$. The authors of \cite{everaert2010homology} introduce the concept of a \emph{protoadditive functor}. A functor $F$ between pointed protomodular categories is protoadditive if it preserves split short exact sequences. Moreover, in \cite{everaert2015protoadditive}, the authors show that a functor between pointed protomodular categories that preserves the zero object is protoadditive if and only if it preserves pullbacks along split epimorphisms. This characterization enables the extension of the notion of protoadditivity to the non-pointed case. 
\begin{definition}
    A functor $S \colon \cat{C} \rightarrow \cat{D}$ between protomodular categories that have both a terminal and an initial object is \emph{protoadditive} if it preserves the terminal object, the initial object, and pullbacks along split epimorphisms.
\end{definition}
\begin{proposition}
The functor $S \colon \mv \rightarrow s\mv$ is protoadditive.
\begin{proof}
It is clear that $S(\terminal)=\terminal/\terminal=\terminal$ and $S(\initial)=\initial/\{0\}=\initial$. Moreover, since $s\mv$ is a full reflective subcategory of $\mv$, we know that $s\mv$ is closed under the formation of limits. We consider, in $\mv$, the pullback of a split epimorphism along an arbitrary morphism
% https://q.uiver.app/?q=WzAsNCxbMCwwLCJBIFxcdGltZXNfQiBDIl0sWzEsMCwiQyJdLFswLDEsIkEiXSxbMSwxLCJCIl0sWzIsMywicCIsMCx7Im9mZnNldCI6LTF9XSxbMywyLCJzIiwwLHsib2Zmc2V0IjotMX1dLFsxLDMsImYiXSxbMCwxLCJcXHBpX0MiXSxbMCwyLCJcXHBpX0EiLDJdLFswLDMsIiIsMix7InN0eWxlIjp7Im5hbWUiOiJjb3JuZXIifX1dXQ==
\[\begin{tikzcd}
	{A \times_B C} & C \\
	A & B.
	\arrow["p", shift left=1, from=2-1, to=2-2]
	\arrow["s", shift left=1, from=2-2, to=2-1]
	\arrow["g", from=1-2, to=2-2]
	\arrow["{\pi_C}", from=1-1, to=1-2]
	\arrow["{\pi_A}"', from=1-1, to=2-1]
	\arrow["\lrcorner"{anchor=center, pos=0.125}, draw=none, from=1-1, to=2-2]
\end{tikzcd}\]
We observe that 
\begin{align*}
    \Rad(A \times_B C)&= \{ (a,c) \in A \times C \, | \, a \in \Rad(A), \, c \in \Rad(C), \, p(a)=g(c) \}\\ 
    &=\Rad(A \times C) \cap (A \times_B C).
\end{align*}
We now proceed to compute the pullback of $S(p)$ along $S(g)$ in $s\mv$:
% https://q.uiver.app/?q=WzAsNSxbMSwxLCJTKEEpIFxcdGltZXNfe1MoQil9IFMoQykiXSxbMiwxLCJTKEMpIl0sWzEsMiwiUyhBKSJdLFsyLDIsIlMoQikiXSxbMCwwLCJTKEEgXFx0aW1lc19CIEMpIl0sWzIsMywiUyhwKSIsMCx7Im9mZnNldCI6LTF9XSxbMywyLCJTKHMpIiwwLHsib2Zmc2V0IjotMX1dLFswLDIsIlxccGlfe1MoQSl9IiwyLHsib2Zmc2V0IjoxfV0sWzEsMywiUyhnKSIsMCx7Im9mZnNldCI6MX1dLFswLDEsIlxccGlfe1MoQyl9IiwwLHsib2Zmc2V0IjotMX1dLFswLDMsIiIsMCx7InN0eWxlIjp7Im5hbWUiOiJjb3JuZXIifX1dLFs0LDIsIlMoXFxwaV9BKSIsMix7ImN1cnZlIjo0fV0sWzQsMSwiUyhcXHBpX0MpIiwwLHsiY3VydmUiOi00fV0sWzQsMCwiXFxleGlzdHMhIFxcdmFycGhpIl1d
\[\begin{tikzcd}
	{S(A \times_B C)} \\
	& {S(A) \times_{S(B)} S(C)} & {S(C)} \\
	& {S(A)} & {S(B).}
	\arrow["{S(p)}", shift left=1, from=3-2, to=3-3]
	\arrow["{S(s)}", shift left=1, from=3-3, to=3-2]
	\arrow["{\pi_{S(A)}}"', shift right=1, from=2-2, to=3-2]
	\arrow["{S(g)}", shift right=1, from=2-3, to=3-3]
	\arrow["{\pi_{S(C)}}", shift left=1, from=2-2, to=2-3]
	\arrow["\lrcorner"{anchor=center, pos=0.125}, draw=none, from=2-2, to=3-3]
	\arrow["{S(\pi_A)}"', curve={height=24pt}, from=1-1, to=3-2]
	\arrow["{S(\pi_C)}", curve={height=-24pt}, from=1-1, to=2-3]
	\arrow["{\exists! \varphi}", from=1-1, to=2-2]
\end{tikzcd}\]
Using the fact that $S(p)S(\pi_A)=S(p\pi_A)=S(g\pi_C)=S(g)S(\pi_C)$, the universal property of this pullback defines an arrow $\varphi \colon S(A \times_B C) \rightarrow S(A) \times_{S(B)} S(C)$ such that $\varphi([a,c])=([a],[c])$. We will prove that $\varphi$ is an isomorphism. It is not difficult to see that $\varphi$ is injective: if $\varphi([a,c])=([0],[0])$, then we have $(a,c) \in A \times_B C$, $a \in \Rad(A)$, and $c \in \Rad(C)$. Thus, $(a,c) \in \Rad(A \times C) \cap (A \times_B C)$, and so $[a,c]=[0,0]$. Now, let us show that $\varphi$ is surjective. Let $([a],[c]) \in S(A) \times_{S(B)} S(C)$. Since $S(p)([a])=S(g)([c])$, we have $[p(a)]=[g(c)]$, and so $[sp(a)]=[sg(c)]$. Our goal is to find $a' \in A$ and $c' \in C$ such that $p(a')=g(c')$ and $[a']=[a]$, $[c']=[c]$. First, we observe that $(sp(a) \ominus a, 0) \in A \times_B C$ and $(a \ominus sp(a), 0) \in A \times_B C$. By the definition of $\varphi$, we get $\varphi([sp(a) \ominus a, 0])=([sp(a) \ominus a],[0])$ and $\varphi([a \ominus sp(a), 0])=([a \ominus sp(a)],[0])$. Additionally, we have $(sg(c),c) \in A \times_B C$ and $\varphi([sg(c),c])=([sg(c)],[c])=([sp(a)],[c])$ (since $[sp(a)]=[sg(c)]$). Let us define $x \coloneqq [a \ominus sp(a), 0]$, $y \coloneqq [sp(a) \ominus a, 0]$,  and $z \coloneqq [sg(c),c]$, for simplicity. Recalling that $$(a \ominus sp(a)) \oplus ((a \oplus \lnot sp(a)) \odot sp(a))=a,$$ we get $$\varphi (x \oplus (\lnot y \odot z))=([a],[c]).$$ Thus, $\varphi$ is surjective, and since we have already shown that it is injective, we conclude that $\varphi$ is an isomorphism.
\end{proof}
\end{proposition}
Additionally, we will prove how the pretorsion theory in $\mv$ studied in this section induces a stable factorization system.
\begin{lemma}\label{nuovo lemma}
Given a morphism $f \colon A \rightarrow B$ in $\mv$ define $$\myfunc{\overline{f}}{A/ \theta_f}{B}{[a]}{f(a),}$$ where $\theta_f \coloneqq \kerr(f) \cap \Rad(A)$ and $\overline{f}$ is well defined since $\theta_f \subseteq \Rad(f)$. Then we have $$\kerr(\overline{f}) \cap \Rad(A/ \theta_f)= \{ 0 \}.$$
\begin{proof}
Fix an element $[a] \in \kerr(\overline{f})  \cap \Rad(A/ \theta_f)$. Since $f(a)=\overline{f}([a])=0$, we immediately get $a \in \kerr(f)$. We consider the quotient projection $\pi_f \colon A \twoheadrightarrow A/ \theta_f$. Given that $[a]=\pi_f(a) \in \Rad(A/ \theta_f)=\bigcap \{ M \subseteq A/\theta_f \, | \, M \text{ is a maximal ideal} \}$ we deduce that $a$ belongs to $$\bigcap \{ \inv{\pi_f}(M) \subseteq A \, | \, M \text{ is a maximal ideal of } A/\theta_f \}.$$ Thanks to \cite{cignoli2013algebraic}, Proposition 1.2.10, we know that $\inv{\pi_f}$ defines a bijection, which preserves and reflects the order, between $\{I \subseteq A/ \theta_f \, | \, I \text{ is an ideal} \}$ and $\{J \subseteq A\, | \, J \text{ is an ideal and } J \supseteq \theta_f \}$. We can observe that the maximal elements of the poset of ideals in $A$ are precisely the maximal elements of the poset $\{J \subseteq A \,|\, J \text{ is an ideal and } J \supseteq \theta_f \}$, as $\theta_f \subseteq \Rad(A)$ and the radical is contained in every maximal ideal. Thus, we have 
\begin{align*}
    &\bigcap \{ \inv{\pi_f}(M) \subseteq A \, | \, M \text{ is a maximal ideal of } A/\theta_f \}=\\ &\bigcap \{ N \subseteq A \, | \, N \text{ is a maximal ideal of } A \text{ and } N \supseteq \theta_f \} = \Rad(A). 
\end{align*}
Consequently, if $[a] \in \kerr(\overline{f}) \cap \Rad(A/ \theta_f)$, then $a \in \kerr(f) \cap \Rad(A)=\theta_f$, and therefore $[a]=0$.
\end{proof}
\end{lemma}

\begin{proposition}
We define the two following classes of arrows in $\mv$ 
\begin{align*}
    \facte &\coloneqq \{ e \colon A \rightarrow B \in \Arr(\mv) \, | \, e \emph{ is surjective, } \kerr(e) \subseteq \Rad(A) \} \textit{ and}\\
    \factm &\coloneqq \{ m \colon A \rightarrow B \in \Arr(\mv) \, | \, \kerr(m) \cap \Rad(A)= \{0 \} \}.
\end{align*}
Then the pair $(\facte, \factm)$ forms a stable factorization system for $\mv$.
\begin{proof}
We start by considering a commutative square
\[\begin{tikzcd}
	A & B \\
	C & D,
	\arrow["{e \in \facte}", from=1-1, to=1-2]
	\arrow["{m \in \factm}"', from=2-1, to=2-2]
	\arrow["g"', from=1-1, to=2-1]
	\arrow["h", from=1-2, to=2-2]
\end{tikzcd}\]
and given that $e$ is surjective, we assume that $B$ is equal to $A / \kerr(e)$ and $e$ is the quotient projection. We observe that for every element $a \in \kerr(e) \subseteq \Rad(A)$, we have $g(a) \in \Rad(C)$ and, furthermore, $mg(a)=he(a)=h(0)=0$. Therefore, we get $g(a) \in \Rad(C) \cap \kerr(m)= \{ 0 \}$ which implies that the arrow $d \colon B \rightarrow C$, where $d([a]) \coloneqq g(a)$, is well defined. Additionally, we can see that $md([a])=mg(a)=he(a)=h([a])$. Finally, $d$ is unique since $e$ is an epimorphism. Next, we consider an arbitrary arrow $f \colon A \rightarrow B$ in $\mv$ and we construct the factorization
% https://q.uiver.app/?q=WzAsMyxbMCwwLCJBIl0sWzIsMCwiQiJdLFsxLDEsIkEvXFxrZXJyKGYpXFxjYXAgXFxSYWQoQSkiXSxbMCwxLCJmIl0sWzAsMiwicSIsMix7ImN1cnZlIjoyLCJzdHlsZSI6eyJoZWFkIjp7Im5hbWUiOiJlcGkifX19XSxbMiwxLCJpIiwyLHsiY3VydmUiOjJ9XV0=
\[\begin{tikzcd}
	A && B \\
	& {A/\kerr(f)\cap \Rad(A),}
	\arrow["f", from=1-1, to=1-3]
	\arrow["q"', two heads, from=1-1, to=2-2]
	\arrow["i"', from=2-2, to=1-3]
\end{tikzcd}\]
where $q$ is the quotient projection and $i([a]) \coloneqq f(a)$. Clearly, since $\kerr(q)= \kerr(f) \cap \Rad(A) \subseteq \Rad(A)$, we have $q \in \facte$. Furthermore, thanks to Lemma \ref{nuovo lemma}, we get $$\kerr(i) \cap \Rad(A/\kerr(f) \cap \Rad(A))= \{ 0 \},$$ and so $i \in \factm$. To conclude, we need to show that the factorization system is stable. For this, we consider the pullback
% https://q.uiver.app/?q=WzAsNCxbMCwwLCJBIFxcdGltZXNfSUMiXSxbMSwwLCJDIl0sWzAsMSwiQSJdLFsxLDEsIkIiXSxbMCwxLCJcXHBpX0MiXSxbMiwzLCJlIFxcaW4gXFxmYWN0ZSIsMix7InN0eWxlIjp7ImhlYWQiOnsibmFtZSI6ImVwaSJ9fX1dLFswLDIsIlxccGlfQSIsMl0sWzEsMywiZyJdLFswLDMsIiIsMSx7InN0eWxlIjp7Im5hbWUiOiJjb3JuZXIifX1dXQ==
\[\begin{tikzcd}
	{A \times_B C} & C \\
	A & B
	\arrow["{\pi_C}", two heads, from=1-1, to=1-2]
	\arrow["{e \in \facte}"', two heads, from=2-1, to=2-2]
	\arrow["{\pi_A}"', from=1-1, to=2-1]
	\arrow["g", from=1-2, to=2-2]
	\arrow["\lrcorner"{anchor=center, pos=0.125}, draw=none, from=1-1, to=2-2]
\end{tikzcd}\]
and we observe that $\pi_C$ is surjective, since it is the pullback of a surjective map. Moreover, $\kerr(\pi_C)= \{ (a,0) \in A \times C \, | \, e(a)=0 \}$, and so every element $(a,0) \in \kerr(\pi_C)$ satisfies $a \in \Rad(A)$. Hence we have $\kerr(\pi_C) \subseteq \Rad(A \times_I C)$.
\end{proof}
\end{proposition}

\section{A Higher Galois Theory for MV-Algebras}\label{Higher Galois Theory for MV-Algebras}

In the final part of this paper, we will explore the higher-order Galois theory determined by the adjunction $S \dashv i$. Moreover, we will apply this analysis to define and examine the commutator of ideal subalgebras in relation with this Galois theory. For an MV-algebra $A \neq \terminal$, a subalgebra $S \subseteq A$ is said to be \emph{ideal} if there is an ideal $I$ of $A$ such that $S= I \cup \lnot I$. If $A= \terminal$, the only subalgebra is $\terminal$ itself, and we assume it is also ideal. To express the commutator of ideal subalgebras, we only need to consider regular epimorphisms (i.e.\ surjective maps) of $\mv$ and $s \mv$, as we will see. In particular, to investigate commutators between ideal subalgebras, we need to use the structure $(S,i, |\Ext \mv|, |\Ext s\mv|)$, where $|\Ext s\mv|$ denotes the class of surjective maps of $s\mv$, and $|\Ext \mv|$ denotes the class of surjective maps of $\mv$. Let us begin by observing that this structure satisfies the necessary conditions to be considered an admissible structure in Galois theory. We know that $s \mv$ is a full and reflective subcategory of $\mv$, and thus closed under limits. Therefore, the set of surjective maps of $\mv$ and $s \mv$ are admissible classes of arrows. Additionally, it is evident that for every surjective map $f$ of $\mv$, $S(f)$ is also surjective. Finally, since every component of the unit of the adjunction $S \dashv i$ is surjective (and the composition of surjective maps is surjective), we can conclude that the Galois structure $(S,i, |\Ext \mv|, |\Ext s\mv|)$ is admissible.
\begin{definition}[\cite{bourn2003denormalized}]

Let $\cat{C}$ be a category with pullbacks. A commutative diagram of regular epimorphism 
% https://q.uiver.app/?q=WzAsNCxbMCwwLCJBIl0sWzAsMSwiQiJdLFsxLDAsIkMiXSxbMSwxLCJEIl0sWzEsMywiYiIsMix7InN0eWxlIjp7ImhlYWQiOnsibmFtZSI6ImVwaSJ9fX1dLFsyLDMsImMiLDAseyJzdHlsZSI6eyJoZWFkIjp7Im5hbWUiOiJlcGkifX19XSxbMCwxLCJmIiwyLHsic3R5bGUiOnsiaGVhZCI6eyJuYW1lIjoiZXBpIn19fV0sWzAsMiwiZyIsMCx7InN0eWxlIjp7ImhlYWQiOnsibmFtZSI6ImVwaSJ9fX1dXQ==
\[\begin{tikzcd}
	A & C \\
	B & D
	\arrow["h"', two heads, from=2-1, to=2-2]
	\arrow["k", two heads, from=1-2, to=2-2]
	\arrow["f"', two heads, from=1-1, to=2-1]
	\arrow["g", two heads, from=1-1, to=1-2]
\end{tikzcd}\]
is a regular pushout if the morphism $\langle f, g \rangle \colon A \rightarrow B \times_D C$ defined by the universal property of the pullback 
% https://q.uiver.app/?q=WzAsNSxbMSwyLCJCIl0sWzIsMiwiRCJdLFsyLDEsIkMiXSxbMSwxLCJCIFxcdGltZXNfRCBDIl0sWzAsMCwiQSJdLFswLDEsImIiLDIseyJzdHlsZSI6eyJoZWFkIjp7Im5hbWUiOiJlcGkifX19XSxbMiwxLCJjIiwwLHsic3R5bGUiOnsiaGVhZCI6eyJuYW1lIjoiZXBpIn19fV0sWzMsMCwiXFxwaV9CIiwyLHsic3R5bGUiOnsiaGVhZCI6eyJuYW1lIjoiZXBpIn19fV0sWzMsMiwiXFxwaV9DIiwwLHsic3R5bGUiOnsiaGVhZCI6eyJuYW1lIjoiZXBpIn19fV0sWzQsMCwiZiIsMix7InN0eWxlIjp7ImhlYWQiOnsibmFtZSI6ImVwaSJ9fX1dLFs0LDIsImciLDAseyJzdHlsZSI6eyJoZWFkIjp7Im5hbWUiOiJlcGkifX19XSxbMywxLCIiLDEseyJzdHlsZSI6eyJuYW1lIjoiY29ybmVyIn19XSxbNCwzLCJcXGxhbmdsZSBmLGcgXFxyYW5nbGUiLDEseyJzdHlsZSI6eyJoZWFkIjp7Im5hbWUiOiJlcGkifX19XV0=
\[\begin{tikzcd}
	A \\
	& {B \times_D C} & C \\
	& B & D
	\arrow["h"', two heads, from=3-2, to=3-3]
	\arrow["k", two heads, from=2-3, to=3-3]
	\arrow["{\pi_B}"', two heads, from=2-2, to=3-2]
	\arrow["{\pi_C}", two heads, from=2-2, to=2-3]
	\arrow["f"', two heads, bend right, from=1-1, to=3-2]
	\arrow["g", two heads, bend left, from=1-1, to=2-3]
	\arrow["\lrcorner"{anchor=center, pos=0.125}, draw=none, from=2-2, to=3-3]
	\arrow["{\langle f,g \rangle}"{description}, two heads, from=1-1, to=2-2]
\end{tikzcd}\]
is a regular epimorphism.
\end{definition}

Applying Lemma \ref{Lemma dei pullback} we immediately get the following characterization:

\begin{proposition}\label{mi serve per la fattorizzazione}

A diagram of regular epimorphisms 
\[\begin{tikzcd}
	A & B \\
	C & D
	\arrow["g"', two heads, from=2-1, to=2-2]
	\arrow["k", two heads, from=1-2, to=2-2]
	\arrow["h"', two heads, from=1-1, to=2-1]
	\arrow["f", two heads, from=1-1, to=1-2]
\end{tikzcd}\]
in the category $\mv$ is a regular pushout if and only if $h(\kerr(f))=\kerr(g)$.
\end{proposition}
Let $\Extmv$ be the category whose objects are the extensions (i.e.\ regular epimorphisms) of MV-algebras and whose morphisms are the commutative diagrams between them; let $\CExtmv$ be the full subcategory of $\Extmv$ determined by the central extension.
We define the functor $S_1$ as follows:
% https://q.uiver.app/?q=WzAsMTAsWzAsMCwiXFxFeHRtdiJdLFs0LDAsIlxcQ0V4dG12Il0sWzAsMSwiQSJdLFszLDEsIkEvXFx0aGV0YV9mIl0sWzMsNCwiQy9cXHRoZXRhX2ciXSxbMSwxLCJCIl0sWzAsNCwiQyJdLFsxLDQsIkQiXSxbNCwxLCJCIl0sWzQsNCwiRCJdLFszLDQsIlxcYmFye2h9IiwyXSxbMCwxLCJTXzEiXSxbMiw1LCJmIiwwLHsic3R5bGUiOnsiaGVhZCI6eyJuYW1lIjoiZXBpIn19fV0sWzYsNywiZyIsMix7InN0eWxlIjp7ImhlYWQiOnsibmFtZSI6ImVwaSJ9fX1dLFsyLDYsImgiLDJdLFs1LDcsImsiXSxbMyw4LCJcXGJhcntmfSIsMCx7InN0eWxlIjp7ImhlYWQiOnsibmFtZSI6ImVwaSJ9fX1dLFs0LDksIlxcYmFye2d9IiwyLHsic3R5bGUiOnsiaGVhZCI6eyJuYW1lIjoiZXBpIn19fV0sWzgsOSwiayJdLFsxNSwxMCwiIiwwLHsic2hvcnRlbiI6eyJzb3VyY2UiOjIwLCJ0YXJnZXQiOjIwfSwic3R5bGUiOnsidGFpbCI6eyJuYW1lIjoibWFwcyB0byJ9fX1dXQ==
\[\begin{tikzcd}[row sep=8pt, column sep=12pt]
	\Extmv &&&& \CExtmv \\
	A & B && {A/\theta_f} & B \\
	\\
	\\
	C & D && {C/\theta_g} & D,
	\arrow[""{name=0, anchor=center, inner sep=0}, "{\overline{h}}"', from=2-4, to=5-4]
	\arrow["{S_1}", from=1-1, to=1-5]
	\arrow["f", two heads, from=2-1, to=2-2]
	\arrow["g"', two heads, from=5-1, to=5-2]
	\arrow["h"', from=2-1, to=5-1]
	\arrow[""{name=1, anchor=center, inner sep=0}, "k", from=2-2, to=5-2]
	\arrow["{\overline{f}}", two heads, from=2-4, to=2-5]
	\arrow["{\overline{g}}"', two heads, from=5-4, to=5-5]
	\arrow["k", from=2-5, to=5-5]
	\arrow[shorten <=13pt, shorten >=13pt, maps to, from=1, to=0]
\end{tikzcd}\]
where $\theta_f \coloneqq \kerr(f) \cap \Rad(A)$ and $\theta_g \coloneqq \kerr(g) \cap \Rad(C)$. Moreover, we define $\overline{f}([a]) \coloneqq f(a)$, $\overline{g}([c]) \coloneqq g(c)$ and $\overline{h}([a]) \coloneqq [h(a)]$. Since $\theta_f \subseteq \kerr(f)$ we immediately get that $\overline{f}$ is well defined (in a similar way one can show that also $\overline{g}$ is well defined). Let us consider an element $a \in \theta_f=\kerr(f) \cap \Rad(A)$; then $h(a) \in \theta_g=\kerr(g) \cap \Rad(C)$ (from $a \in \Rad(A)$ we get $h(a) \in \Rad(C)$ and from $a \in \kerr(f)$ we obtain $h(a) \in \kerr(g)$ because $gh=kf$); therefore, $\overline{h}$ is well defined. It remains to prove that $\overline{f} \colon A/ \theta_f \twoheadrightarrow B$ and $\overline{g} \colon C/ \theta_g \twoheadrightarrow D$ are central extensions. In other words, we have to show that $$\kerr(\overline{f}) \cap \Rad(A/ \theta_f)= \{ 0 \} \text{ and }\kerr(\overline{g}) \cap \Rad(C/ \theta_g)= \{ 0 \}.$$
These equalities are a direct consequence of Lemma \ref{nuovo lemma}.\\

Being $\mv$ a Mal'tsev variety, it was proved in \cite{gran2004galois} that the inclusion functor $i_1 \colon \CExtmv \rightarrow \Extmv$ has a left adjoint, which is, in our case, $S_1$. This adjunction gives rise to the admissible Galois structure $$\Gamma_1 \coloneqq (S_1, i_1, \Extdmv, \ExtCExtmv);$$where the class of arrows $\Extdmv$ of $\Extmv$ is the class of squares which are regular pushouts, and the class of arrows $\ExtCExtmv$ of $\CExtmv$ is given by the squares of $\Extdmv$ for which the horizontal arrows are central extensions.\\

We propose a study of the central extensions for the structure $\Gamma_1$.

\begin{proposition}

Consider an element $(h,k) \in \Extdmv$
% https://q.uiver.app/?q=WzAsNCxbMCwwLCJBIl0sWzEsMCwiQiJdLFswLDEsIkMiXSxbMSwxLCJEIl0sWzIsMywiZyIsMix7InN0eWxlIjp7ImhlYWQiOnsibmFtZSI6ImVwaSJ9fX1dLFsxLDMsImsiLDAseyJzdHlsZSI6eyJoZWFkIjp7Im5hbWUiOiJlcGkifX19XSxbMCwxLCJmIiwwLHsic3R5bGUiOnsiaGVhZCI6eyJuYW1lIjoiZXBpIn19fV0sWzAsMiwiaCIsMix7InN0eWxlIjp7ImhlYWQiOnsibmFtZSI6ImVwaSJ9fX1dXQ==
\[\begin{tikzcd}
	A & B \\
	C & D;
	\arrow["g"', two heads, from=2-1, to=2-2]
	\arrow["k", two heads, from=1-2, to=2-2]
	\arrow["f", two heads, from=1-1, to=1-2]
	\arrow["h"', two heads, from=1-1, to=2-1]
\end{tikzcd}\]
$(h,k)$ is a normal extension for $\Gamma_1$ if and only if $\langle h,f \rangle \colon A \twoheadrightarrow C \times_D B$ is a normal (or central) extension for $\Gamma$.

\begin{proof}

Let us consider the pullback of $(h,k)$ along $(h,k)$: 
% https://q.uiver.app/?q=WzAsOCxbMCwzLCJBIl0sWzIsMywiQyJdLFsxLDIsIkIiXSxbMywyLCJEIl0sWzAsMSwiQSBcXHRpbWVzX0MgQSJdLFsyLDEsIkEiXSxbMSwwLCJCIFxcdGltZXNfRCBCIl0sWzMsMCwiQiJdLFswLDIsImYiLDIseyJzdHlsZSI6eyJoZWFkIjp7Im5hbWUiOiJlcGkifX19XSxbMCwxLCJoIiwyLHsic3R5bGUiOnsiaGVhZCI6eyJuYW1lIjoiZXBpIn19fV0sWzEsMywiZyIsMix7InN0eWxlIjp7ImhlYWQiOnsibmFtZSI6ImVwaSJ9fX1dLFs0LDUsIlxccGlfQV4yIiwyLHsibGFiZWxfcG9zaXRpb24iOjIwLCJzdHlsZSI6eyJoZWFkIjp7Im5hbWUiOiJlcGkifX19XSxbNCwwLCJcXHBpX0FeMSIsMl0sWzUsMSwiaCIsMix7ImxhYmVsX3Bvc2l0aW9uIjozMH1dLFs0LDYsImYgXFx0aW1lcyBmIiwwLHsibGFiZWxfcG9zaXRpb24iOjMwLCJzdHlsZSI6eyJoZWFkIjp7Im5hbWUiOiJlcGkifX19XSxbNSw3LCJmIiwyLHsic3R5bGUiOnsiaGVhZCI6eyJuYW1lIjoiZXBpIn19fV0sWzcsMywiayIsMCx7InN0eWxlIjp7ImhlYWQiOnsibmFtZSI6ImVwaSJ9fX1dLFs2LDIsIlxccGlfQl4xIiwwLHsibGFiZWxfcG9zaXRpb24iOjMwLCJzdHlsZSI6eyJoZWFkIjp7Im5hbWUiOiJlcGkifX19XSxbMiwzLCJrIiwyLHsibGFiZWxfcG9zaXRpb24iOjIwLCJzdHlsZSI6eyJoZWFkIjp7Im5hbWUiOiJlcGkifX19XSxbNiw3LCJcXHBpX0FeMSIsMCx7InN0eWxlIjp7ImhlYWQiOnsibmFtZSI6ImVwaSJ9fX1dLFs0LDEsIiIsMix7InN0eWxlIjp7Im5hbWUiOiJjb3JuZXIifX1dLFs2LDMsIiIsMSx7InN0eWxlIjp7Im5hbWUiOiJjb3JuZXIifX1dXQ==
\[\begin{tikzcd}
	& {B \times_D B} && B \\
	{A \times_C A} && A \\
	& B && D \\
	A && C.
	\arrow["f"', two heads, from=4-1, to=3-2]
	\arrow["h"', two heads, from=4-1, to=4-3]
	\arrow["g"', two heads, from=4-3, to=3-4]
	\arrow["{\pi_A^2}"'{pos=0.2}, two heads, from=2-1, to=2-3]
	\arrow["{\pi_A^1}"', from=2-1, to=4-1]
	\arrow["h"', two heads, {pos=0.3}, from=2-3, to=4-3]
	\arrow["{f \times f}"{pos=0.3}, two heads, from=2-1, to=1-2]
	\arrow["f"', two heads, from=2-3, to=1-4]
	\arrow["k", two heads, from=1-4, to=3-4]
	\arrow["{\pi_B^1}"{pos=0.3}, two heads, from=1-2, to=3-2]
	\arrow["k"'{pos=0.2}, two heads, from=3-2, to=3-4]
	\arrow["{\pi_B^1}", two heads, from=1-2, to=1-4]
	\arrow["\lrcorner"{anchor=center, pos=0.125}, draw=none, from=2-1, to=4-3]
	\arrow["\lrcorner"{anchor=center, pos=0.125}, draw=none, from=1-2, to=3-4]
\end{tikzcd}\]
Given an object $f \colon A \twoheadrightarrow B$ of $\Extmv$, we observe that the $f$-component $\eta^1_f$ of the unit of the adjunction $S_1 \dashv i_1$ is given by
% https://q.uiver.app/?q=WzAsNCxbMCwwLCJBIl0sWzEsMCwiQiJdLFswLDEsIkEvXFx0aGV0YV9mIl0sWzEsMSwiQi4iXSxbMCwxLCJmIiwwLHsic3R5bGUiOnsiaGVhZCI6eyJuYW1lIjoiZXBpIn19fV0sWzAsMiwiXFxwaV9mIiwyLHsic3R5bGUiOnsiaGVhZCI6eyJuYW1lIjoiZXBpIn19fV0sWzIsMywiXFxiYXJ7Zn0iLDJdXQ==
\[\begin{tikzcd}
	A & B \\
	{A/\theta_f} & {B.}
	\arrow["f", two heads, from=1-1, to=1-2]
	\arrow["{\pi_f}"', two heads, from=1-1, to=2-1]
	\arrow["{\overline{f}}"', two heads, from=2-1, to=2-2]
	\arrow[equal, from=1-2, to=2-2]
\end{tikzcd}\]
In other words, $\eta^1_f=(\pi_f, id_B)$.
We construct the diagram associated with the naturality of $\eta^1$
% https://q.uiver.app/?q=WzAsOCxbMCwzLCJBIl0sWzIsMywiQS8gXFx0aGV0YV9mIl0sWzEsMiwiQiJdLFszLDIsIkIiXSxbMCwxLCJBIFxcdGltZXNfQyBBIl0sWzIsMSwiQSBcXHRpbWVzX0MgQS8gXFx0aGV0YV97ZiBcXHRpbWVzIGZ9Il0sWzEsMCwiQiBcXHRpbWVzX0QgQiJdLFszLDAsIkIgXFx0aW1lc19EIEIiXSxbMCwyLCJmIiwyLHsic3R5bGUiOnsiaGVhZCI6eyJuYW1lIjoiZXBpIn19fV0sWzAsMSwiXFxwaV9mIiwyLHsic3R5bGUiOnsiaGVhZCI6eyJuYW1lIjoiZXBpIn19fV0sWzEsMywiXFxvdmVybGluZXtmfSIsMix7InN0eWxlIjp7ImhlYWQiOnsibmFtZSI6ImVwaSJ9fX1dLFs0LDUsIlxccGlfe2YgXFx0aW1lcyBmfSIsMix7ImxhYmVsX3Bvc2l0aW9uIjoyMCwic3R5bGUiOnsiaGVhZCI6eyJuYW1lIjoiZXBpIn19fV0sWzQsMCwiXFxwaV9BXjEiLDJdLFs1LDEsImgiLDIseyJsYWJlbF9wb3NpdGlvbiI6MzB9XSxbNCw2LCJmIFxcdGltZXMgZiIsMCx7ImxhYmVsX3Bvc2l0aW9uIjozMCwic3R5bGUiOnsiaGVhZCI6eyJuYW1lIjoiZXBpIn19fV0sWzUsNywiXFxvdmVybGluZXtmIFxcdGltZXMgZn0iLDIseyJzdHlsZSI6eyJoZWFkIjp7Im5hbWUiOiJlcGkifX19XSxbNywzLCJrIiwwLHsic3R5bGUiOnsiaGVhZCI6eyJuYW1lIjoiZXBpIn19fV0sWzYsMiwiXFxwaV9CXjEiLDAseyJsYWJlbF9wb3NpdGlvbiI6MzAsInN0eWxlIjp7ImhlYWQiOnsibmFtZSI6ImVwaSJ9fX1dLFsyLDMsIiIsMSx7InN0eWxlIjp7ImhlYWQiOnsibmFtZSI6Im5vbmUifX19XSxbNiw3LCIiLDEseyJzdHlsZSI6eyJoZWFkIjp7Im5hbWUiOiJub25lIn19fV1d
\[\begin{tikzcd}
	& {B \times_D B} && {B \times_D B} \\
	{A \times_C A} && {A \times_C A/ \theta_{f \times f}} \\
	& B && B \\
	A && {A/ \theta_f.}
	\arrow[equal, from=3-2, to=3-4]
	\arrow[equal, from=1-2, to=1-4]
	\arrow["f"', two heads, from=4-1, to=3-2]
	\arrow["{\pi_f}"', two heads, from=4-1, to=4-3]
	\arrow["{\overline{f}}"', two heads, from=4-3, to=3-4]
	\arrow["{\pi_{f \times f}}"'{pos=0.2}, two heads, from=2-1, to=2-3]
	\arrow["{\pi_A^1}"', from=2-1, to=4-1]
	\arrow["\overline{\pi_A^1}"'{pos=0.3}, from=2-3, to=4-3]
	\arrow["{f \times f}"{pos=0.3}, two heads, from=2-1, to=1-2]
	\arrow["{\overline{f \times f}}"', two heads, from=2-3, to=1-4]
	\arrow["k", two heads, from=1-4, to=3-4]
	\arrow["{\pi_B^1}"{pos=0.3}, two heads, from=1-2, to=3-2]
\end{tikzcd}\]
We have to prove that this cube is a pullback (i.e.\ the front square is a pullback) if and only if $\langle h, f \rangle$ is central for $\Gamma$ (i.e.\ $\kerr(\langle h, f \rangle) \cap \Rad(A)= \{ 0 \}$).\\
We know that the square
% https://q.uiver.app/?q=WzAsNCxbMCwxLCJBIl0sWzEsMSwiQS8gXFx0aGV0YV9mIl0sWzAsMCwiQSBcXHRpbWVzX0MgQSJdLFsxLDAsIkEgXFx0aW1lc19DIEEvIFxcdGhldGFfe2YgXFx0aW1lcyBmfSJdLFswLDEsIlxccGlfZiIsMix7InN0eWxlIjp7ImhlYWQiOnsibmFtZSI6ImVwaSJ9fX1dLFsyLDMsIlxccGlfe2YgXFx0aW1lcyBmfSIsMCx7ImxhYmVsX3Bvc2l0aW9uIjoyMCwic3R5bGUiOnsiaGVhZCI6eyJuYW1lIjoiZXBpIn19fV0sWzIsMCwiXFxwaV9BXjEiLDJdLFszLDEsImgiLDAseyJsYWJlbF9wb3NpdGlvbiI6MzB9XV0=
\[\begin{tikzcd}
	{A \times_C A} & {A \times_C A/ \theta_{f \times f}} \\
	A & {A/ \theta_f}
	\arrow["{\pi_f}"', two heads, from=2-1, to=2-2]
	\arrow["{\pi_{f \times f}}", two heads, from=1-1, to=1-2]
	\arrow["{\pi_A^1}"', from=1-1, to=2-1]
	\arrow["\overline{\pi_A^1}", from=1-2, to=2-2]
\end{tikzcd}\]
is a pullback if and only if the restriction $\pi_A^1 \colon \theta_{f \times f} \rightarrow \theta_f$ is a bijection (thanks to Lemma \ref{Lemma dei pullback}). We prove that the restriction of $\pi_A^1$ is always surjective. Let $a \in \theta_f=\kerr(f) \cap \Rad(A)$; then $(a,a) \in \kerr(f \times f) \cap \Rad(A \times A)$, and $\pi_A^1(a,a)=a$. By Lemma \ref{ideali}, we know that the restriction of $\pi_A^1$ is injective if and only if $\kerr(\pi_A^1) \cap \theta_{f \times f}=\{ 0 \}$. Specifically, $\kerr(\pi_A^1) \cap \theta_{f \times f}=\kerr(\pi_A^1) \cap \kerr(f \times f) \cap \Rad(A \times_C A)=\{(0,a) \in A \times A \, | \, f(a)=0, h(a)=0, a \in \Rad(A) \}$. It is clear that $\kerr(\pi_A^1) \cap \theta_{f \times f}=\{ 0 \}$ if and only if $\kerr(f) \cap \kerr(h) \cap \Rad(A)=\{ 0 \}$. The last statement is true if and only if $\langle h,f \rangle \colon A \twoheadrightarrow C \times_D B$ is a normal extension for $\Gamma$, observing that $\kerr(\langle h,f \rangle)=\kerr(h) \cap \kerr(f)$.
\end{proof}

\end{proposition}

\begin{proposition}

Consider an element $(h,k) \in \Extdmv$
% https://q.uiver.app/?q=WzAsNCxbMCwwLCJBIl0sWzEsMCwiQiJdLFswLDEsIkMiXSxbMSwxLCJEIl0sWzIsMywiZyIsMix7InN0eWxlIjp7ImhlYWQiOnsibmFtZSI6ImVwaSJ9fX1dLFsxLDMsImsiLDAseyJzdHlsZSI6eyJoZWFkIjp7Im5hbWUiOiJlcGkifX19XSxbMCwxLCJmIiwwLHsic3R5bGUiOnsiaGVhZCI6eyJuYW1lIjoiZXBpIn19fV0sWzAsMiwiaCIsMix7InN0eWxlIjp7ImhlYWQiOnsibmFtZSI6ImVwaSJ9fX1dXQ==
\[\begin{tikzcd}
	A & B \\
	C & D;
	\arrow["g"', two heads, from=2-1, to=2-2]
	\arrow["k", two heads, from=1-2, to=2-2]
	\arrow["f", two heads, from=1-1, to=1-2]
	\arrow["h"', two heads, from=1-1, to=2-1]
\end{tikzcd}\]
$(h,k)$ is a central extension for $\Gamma_1$ if and only if $(h,k)$ is a normal extension for $\Gamma_1$.
\begin{proof}

Clearly if $(h,k)$ is normal then it is central. Let us prove that also the other implication is true. We consider two objects $(h,k)$ and $(h',k')$ of $\Extdmv$, defined by the following squares:
% https://q.uiver.app/?q=WzAsOCxbMCwwLCJBIl0sWzEsMCwiQiJdLFswLDEsIkMiXSxbMSwxLCJEIl0sWzIsMCwiQSciXSxbMiwxLCJDIl0sWzMsMCwiQiciXSxbMywxLCJEIl0sWzAsMSwiZiIsMCx7InN0eWxlIjp7ImhlYWQiOnsibmFtZSI6ImVwaSJ9fX1dLFswLDIsImgiLDIseyJzdHlsZSI6eyJoZWFkIjp7Im5hbWUiOiJlcGkifX19XSxbMiwzLCJnIiwyLHsic3R5bGUiOnsiaGVhZCI6eyJuYW1lIjoiZXBpIn19fV0sWzEsMywiayIsMCx7InN0eWxlIjp7ImhlYWQiOnsibmFtZSI6ImVwaSJ9fX1dLFs0LDUsImgnIiwyLHsic3R5bGUiOnsiaGVhZCI6eyJuYW1lIjoiZXBpIn19fV0sWzQsNiwiZiciLDAseyJzdHlsZSI6eyJoZWFkIjp7Im5hbWUiOiJlcGkifX19XSxbNiw3LCJrJyIsMCx7InN0eWxlIjp7ImhlYWQiOnsibmFtZSI6ImVwaSJ9fX1dLFs1LDcsImciLDIseyJzdHlsZSI6eyJoZWFkIjp7Im5hbWUiOiJlcGkifX19XV0=
\[\begin{tikzcd}
	A & B & {A'} & {B'} \\
	C & D & C & D.
	\arrow["f", two heads, from=1-1, to=1-2]
	\arrow["h"', two heads, from=1-1, to=2-1]
	\arrow["g"', two heads, from=2-1, to=2-2]
	\arrow["k", two heads, from=1-2, to=2-2]
	\arrow["{h'}"', two heads, from=1-3, to=2-3]
	\arrow["{f'}", two heads, from=1-3, to=1-4]
	\arrow["{k'}", two heads, from=1-4, to=2-4]
	\arrow["g"', two heads, from=2-3, to=2-4]
\end{tikzcd}\]
The pullback of $(h,k)$ along $(h',k')$ is given by the commutative cube
% https://q.uiver.app/?q=WzAsOCxbMCwzLCJBIl0sWzIsMywiQyJdLFsxLDIsIkIiXSxbMywyLCJEIl0sWzAsMSwiQSBcXHRpbWVzX0MgQSciXSxbMiwxLCJBJyJdLFsxLDAsIkIgXFx0aW1lc19EIEInIl0sWzMsMCwiQiciXSxbMCwyLCJmIiwyLHsic3R5bGUiOnsiaGVhZCI6eyJuYW1lIjoiZXBpIn19fV0sWzAsMSwiaCIsMix7InN0eWxlIjp7ImhlYWQiOnsibmFtZSI6ImVwaSJ9fX1dLFsxLDMsImciLDIseyJzdHlsZSI6eyJoZWFkIjp7Im5hbWUiOiJlcGkifX19XSxbNCw1LCJcXHBpX3tBJ30iLDIseyJsYWJlbF9wb3NpdGlvbiI6MjAsInN0eWxlIjp7ImhlYWQiOnsibmFtZSI6ImVwaSJ9fX1dLFs0LDAsIlxccGlfQSIsMl0sWzUsMSwiaCciLDIseyJsYWJlbF9wb3NpdGlvbiI6MzB9XSxbNCw2LCJmIFxcdGltZXMgZiciLDAseyJsYWJlbF9wb3NpdGlvbiI6MzAsInN0eWxlIjp7ImhlYWQiOnsibmFtZSI6ImVwaSJ9fX1dLFs1LDcsImYnIiwyLHsic3R5bGUiOnsiaGVhZCI6eyJuYW1lIjoiZXBpIn19fV0sWzcsMywiayciLDAseyJzdHlsZSI6eyJoZWFkIjp7Im5hbWUiOiJlcGkifX19XSxbNiwyLCJcXHBpX0IiLDAseyJsYWJlbF9wb3NpdGlvbiI6MzAsInN0eWxlIjp7ImhlYWQiOnsibmFtZSI6ImVwaSJ9fX1dLFsyLDMsImsiLDIseyJsYWJlbF9wb3NpdGlvbiI6MjAsInN0eWxlIjp7ImhlYWQiOnsibmFtZSI6ImVwaSJ9fX1dLFs2LDcsIlxccGlfe0InfSIsMCx7InN0eWxlIjp7ImhlYWQiOnsibmFtZSI6ImVwaSJ9fX1dLFs0LDEsIiIsMix7InN0eWxlIjp7Im5hbWUiOiJjb3JuZXIifX1dLFs2LDMsIiIsMSx7InN0eWxlIjp7Im5hbWUiOiJjb3JuZXIifX1dXQ==
\[\begin{tikzcd}
	& {B \times_D B'} && {B'} \\
	{A \times_C A'} && {A'} \\
	& B && D \\
	A && C.
	\arrow["f"', two heads, from=4-1, to=3-2]
	\arrow["h"', two heads, from=4-1, to=4-3]
	\arrow["g"', two heads, from=4-3, to=3-4]
	\arrow["{\pi_{A'}}"'{pos=0.2}, two heads, from=2-1, to=2-3]
	\arrow["{\pi_A}"', from=2-1, to=4-1]
	\arrow["{h'}"'{pos=0.3}, from=2-3, to=4-3]
	\arrow["{f \times f'}"{pos=0.3}, two heads, from=2-1, to=1-2]
	\arrow["{f'}"', two heads, from=2-3, to=1-4]
	\arrow["{k'}", two heads, from=1-4, to=3-4]
	\arrow["{\pi_B}"{pos=0.3}, two heads, from=1-2, to=3-2]
	\arrow["k"'{pos=0.2}, two heads, from=3-2, to=3-4]
	\arrow["{\pi_{B'}}", two heads, from=1-2, to=1-4]
	\arrow["\lrcorner"{anchor=center, pos=0.125}, draw=none, from=2-1, to=4-3]
	\arrow["\lrcorner"{anchor=center, pos=0.125}, draw=none, from=1-2, to=3-4]
\end{tikzcd}\]
We construct the diagram associated with the naturality of $\eta^1$:
% https://q.uiver.app/?q=WzAsOCxbMCwzLCJBJyJdLFsyLDMsIkEvIFxcdGhldGFfe2YnfSJdLFsxLDIsIkInIl0sWzMsMiwiQiciXSxbMCwxLCJBIFxcdGltZXNfQyBBJyJdLFsyLDEsIkEgXFx0aW1lc19DIEEvIFxcdGhldGFfe2YgXFx0aW1lcyBmfSJdLFsxLDAsIkIgXFx0aW1lc19EIEInIl0sWzMsMCwiQiBcXHRpbWVzX0QgQiciXSxbMCwyLCJmJyIsMix7InN0eWxlIjp7ImhlYWQiOnsibmFtZSI6ImVwaSJ9fX1dLFswLDEsIlxccGlfe2YnfSIsMix7InN0eWxlIjp7ImhlYWQiOnsibmFtZSI6ImVwaSJ9fX1dLFsxLDMsIlxcb3ZlcmxpbmV7Zid9IiwyLHsic3R5bGUiOnsiaGVhZCI6eyJuYW1lIjoiZXBpIn19fV0sWzQsNSwiXFxwaV97ZiBcXHRpbWVzIGYnfSIsMix7ImxhYmVsX3Bvc2l0aW9uIjoyMCwic3R5bGUiOnsiaGVhZCI6eyJuYW1lIjoiZXBpIn19fV0sWzQsMCwiXFxwaV97QSd9IiwyXSxbNSwxLCJcXG92ZXJsaW5le1xccGlfe0EnfX0iLDIseyJsYWJlbF9wb3NpdGlvbiI6MzB9XSxbNCw2LCJmIFxcdGltZXMgZiciLDAseyJsYWJlbF9wb3NpdGlvbiI6MzAsInN0eWxlIjp7ImhlYWQiOnsibmFtZSI6ImVwaSJ9fX1dLFs1LDcsIlxcb3ZlcmxpbmV7ZiBcXHRpbWVzIGYnfSIsMix7InN0eWxlIjp7ImhlYWQiOnsibmFtZSI6ImVwaSJ9fX1dLFs3LDMsIlxccGlfe0InfSIsMCx7InN0eWxlIjp7ImhlYWQiOnsibmFtZSI6ImVwaSJ9fX1dLFs2LDIsIlxccGlfe0InfSIsMCx7ImxhYmVsX3Bvc2l0aW9uIjozMCwic3R5bGUiOnsiaGVhZCI6eyJuYW1lIjoiZXBpIn19fV0sWzIsMywiIiwxLHsic3R5bGUiOnsiaGVhZCI6eyJuYW1lIjoibm9uZSJ9fX1dLFs2LDcsIiIsMSx7InN0eWxlIjp7ImhlYWQiOnsibmFtZSI6Im5vbmUifX19XV0=
\[\begin{tikzcd}
	& {B \times_D B'} && {B \times_D B'} \\
	{A \times_C A'} && {A \times_C A/ \theta_{f \times f}} \\
	& {B'} && {B'} \\
	{A'} && {A/ \theta_{f'}.}
	\arrow["{f'}"', two heads, from=4-1, to=3-2]
	\arrow["{\pi_{f'}}"', two heads, from=4-1, to=4-3]
	\arrow["{\overline{f'}}"', two heads, from=4-3, to=3-4]
	\arrow["{\pi_{f \times f'}}"'{pos=0.2}, two heads, from=2-1, to=2-3]
	\arrow["{\pi_{A'}}"', from=2-1, to=4-1]
	\arrow["{\overline{\pi_{A'}}}"'{pos=0.3}, from=2-3, to=4-3]
	\arrow["{f \times f'}"{pos=0.3}, two heads, from=2-1, to=1-2]
	\arrow["{\overline{f \times f'}}"', two heads, from=2-3, to=1-4]
	\arrow["{\pi_{B'}}", two heads, from=1-4, to=3-4]
	\arrow["{\pi_{B'}}"{pos=0.3}, two heads, from=1-2, to=3-2]
	\arrow[equal, from=3-2, to=3-4]
	\arrow[equal, from=1-2, to=1-4]
\end{tikzcd}\]
If the restriction of $\pi_{A'} \colon \theta_{f \times f'} \rightarrow \theta_{f'}$ is bijective, then it is injective and so, thanks to Lemma \ref{ideali}, we obtain $\kerr(\pi_{A'}) \cap \kerr(f \times f') \cap \Rad(A \times_C A')= \{0\}$. But $\kerr(\pi_{A'}) \cap \kerr(f \times f') \cap \Rad(A \times_C A')=\{ (a,0) \in A \times A \, | \, h(a)=0, f(a)=0, a \in \Rad(A) \}$. It is clear that $\kerr(\pi_{A'}) \cap \kerr(f \times f') \cap \Rad(A \times_C A')= \{0\}$ if and only if $\kerr(h) \cap \kerr(f) \cap \Rad(A)= \{0\}$. The last statement is true if and only if $\langle h,f \rangle \colon A \twoheadrightarrow C \times_D B$ is a normal extension for $\Gamma$ (this is implied by the fact that $\kerr(\langle h,f \rangle)=\kerr(h) \cap \kerr(f)$).
\end{proof}

\end{proposition}

\begin{corollary}

Consider an element $(h,k) \in \Extdmv$
% https://q.uiver.app/?q=WzAsNCxbMCwwLCJBIl0sWzEsMCwiQiJdLFswLDEsIkMiXSxbMSwxLCJEIl0sWzIsMywiZyIsMix7InN0eWxlIjp7ImhlYWQiOnsibmFtZSI6ImVwaSJ9fX1dLFsxLDMsImsiLDAseyJzdHlsZSI6eyJoZWFkIjp7Im5hbWUiOiJlcGkifX19XSxbMCwxLCJmIiwwLHsic3R5bGUiOnsiaGVhZCI6eyJuYW1lIjoiZXBpIn19fV0sWzAsMiwiaCIsMix7InN0eWxlIjp7ImhlYWQiOnsibmFtZSI6ImVwaSJ9fX1dXQ==
\[\begin{tikzcd}
	A & B \\
	C & D.
	\arrow["g"', two heads, from=2-1, to=2-2]
	\arrow["k", two heads, from=1-2, to=2-2]
	\arrow["f", two heads, from=1-1, to=1-2]
	\arrow["h"', two heads, from=1-1, to=2-1]
\end{tikzcd}\]
 $(h,k)$ is a central extension for $\Gamma_1$ if and only if $\langle h,f \rangle \colon A \twoheadrightarrow C \times_D B$ is a central extension for $\Gamma$.

\end{corollary}

We are now prepared to describe the commutator of two ideal subalgebras with respect to the Galois structure that has just been studied. To accomplish this, we require the following results.

\begin{lemma}\label{intersezioni banali}

Let $A \neq \terminal$ be an MV-algebra and $I \subseteq A$ a proper ideal. Then $\Rad(A) \cap \lnot I= \emptyset$.

\begin{proof}

Suppose $x \in \Rad(A) \cap \lnot I$. Since $x \in \Rad(A)$, we have $x \leq \lnot x$. However, since $\lnot x \in I$ by assumption, we also have $x \in I$. Thus, we conclude that $1 = x \oplus \lnot x \in I$, which is in contradiction with the fact that $I$ is a proper ideal. Therefore, $\Rad(A) \cap \lnot I$ must be empty, as claimed. 
\end{proof}

\end{lemma}

\begin{lemma}\label{riesz MV}

Let $A$ be an MV-algebra. Consider a triple of elements $x,y,z\in A$ such that $x\leq y\oplus z$. Then, there exist $y_1,z_1\in A$ such that $y_1\leq y$, $z_1\leq z$, and $x=y_1\oplus z_1$.

\begin{proof}

 By \cite{mundici1986interpretation}, Theorem 3.9, we know that $A=[0,u]$ for a certain object $(G,u)$ in $\ulab$, where $\ulab$ denotes the category of lattice-ordered abelian groups with order-unit. We recall that in $A$ the operation $\oplus$ is defined as $a\oplus b=(a+b)\land u$, where $+$ denotes the addition in $(G,u)$. Given $x,y,z\in A$ such that $x\leq y\oplus z \leq y+z$, we can apply the Riesz decomposition property (see, for example, Theorem 2.1 in \cite{librogruppiret}) to obtain two elements $0\leq y_1\leq y$ and $0\leq z_1\leq z$ such that $x=y_1+z_1$. Finally, we observe that $y_1\oplus z_1=(y_1+z_1)\land u=x\land u=x=y_1+z_1$, and thus the statement holds.
\end{proof}
    
\end{lemma}

\begin{proposition}

Let $A$ be an MV-algebra and let $I$ and $J$ be two ideals of $A$. We define $$I \oplus J \coloneqq \{x \in A \, | \, \exists i \in I, \exists j \in J \emph{ such that } x = i \oplus j\}.$$ Then, we have $I \lor J = I \oplus J$, where $I \lor J$ denotes the join of $I$ and $J$ computed in $\Ideals(A)$.

\begin{proof}

To prove that $I \oplus J$ is an ideal of $A$, it suffices to show that $0 \in I \oplus J$, $I \oplus J$ is closed under $\oplus$, and $I \oplus J$ is downward-closed.
Since $0 = 0 \oplus 0$, we have $0 \in I \oplus J$. Moreover, if $x = i_1 \oplus j_1$ and $y = i_2 \oplus j_2$ with $i_1, i_2 \in I$ and $j_1, j_2 \in J$, then $x \oplus y = (i_1 \oplus i_2) \oplus (j_1 \oplus j_2) \in I \oplus J$, which shows that $I \oplus J$ is closed under $\oplus$.
Finally, suppose $z \in I \oplus J$ and $w \in A$ such that $w \leq z$. Then there exist $i \in I$ and $j \in J$ such that $z = i \oplus j$. Applying Lemma \ref{riesz MV}, we obtain $w = i_w \oplus j_w$ with $i_w \leq i$ and $j_w \leq j$. Therefore $w \in I \oplus J$, which shows that $I \oplus J$ is downward-closed. Clearly, any ideal containing both $I$ and $J$ also contains $I \oplus J$, and this completes the proof.
\end{proof}
    
\end{proposition}

Consider an MV-algebra $A$, with $A \neq \terminal$, and fix a pair $I,J$ of proper ideals of $A$. Let $M = I \cup \lnot I$ and $N = J \cup \lnot J$ be the two ideal subalgebras associated with $I$ and $J$, respectively. We aim to show that $$M \lor N= K \cup \lnot K,$$ where $M \lor N$ denotes the subalgebra generated by $M$ and $N$, and $K= I \lor J$ (with this join computed in the poset $\Ideals(A)$). In fact, as $K$ contains both $I$ and $J$, it is clear that $K \cup \lnot K \supseteq M$ and $K \cup \lnot K \supseteq N$, which implies $K \cup \lnot K \supseteq M \lor N$. Moreover, since $M \lor N$ contains both $I$ and $J$, we have $M \lor N \supseteq I \oplus J = I \lor J$ (the last equality holds applying the previous proposition). As $M \lor N$ is an MV-algebra, we also have $M \lor N \supseteq K \cup \lnot K$. Therefore, the equality we need to show holds.\\
In order to describe the commutator between $M$ and $N$, we have to distinguish two cases. If $I \lor J = A$, then the diagram 
% https://q.uiver.app/?q=WzAsNCxbMCwwLCJBIl0sWzEsMCwiQS9KIl0sWzEsMSwiXFx0ZXJtaW5hbCJdLFswLDEsIkEvSSJdLFswLDEsInFfSiIsMCx7InN0eWxlIjp7ImhlYWQiOnsibmFtZSI6ImVwaSJ9fX1dLFsxLDIsIiIsMCx7InN0eWxlIjp7ImhlYWQiOnsibmFtZSI6ImVwaSJ9fX1dLFszLDIsIiIsMix7InN0eWxlIjp7ImhlYWQiOnsibmFtZSI6ImVwaSJ9fX1dLFswLDMsInFfSSIsMix7InN0eWxlIjp7ImhlYWQiOnsibmFtZSI6ImVwaSJ9fX1dXQ==
\[\begin{tikzcd}\label{diagramma IvJ=A}\tag{$\clubsuit$}
	A & {A/J} \\
	{A/I} & \terminal
	\arrow["{q_J}", two heads, from=1-1, to=1-2]
	\arrow[two heads, from=1-2, to=2-2]
	\arrow[two heads, from=2-1, to=2-2]
	\arrow["{q_I}"', two heads, from=1-1, to=2-1]
\end{tikzcd}\]
is a regular pushout. To show this, we only need to prove that the restriction $q_I \colon J \rightarrow A/I$ is surjective. Let $[a] \in A/I$ be an arbitrary element. Since $I \lor J = A$, there exist $i \in I$ and $j \in J$ such that $a = i \oplus j$. Therefore, $[a] = [i \oplus j] = [j] = q_I(j)$. This implies that the restriction of $q_I$ is surjective, completing the proof. If $I \lor J \neq A$, then $K$ is a proper ideal of $A$. Hence, it is easy to see that $K \cap \lnot K = \emptyset$. Therefore, we define a morphism of MV-algebras $\chi \colon M \lor N=K \cup \lnot K \rightarrow \initial$, setting $\chi(x)=1$ if and only if $x \in \lnot K$. Since $I \subseteq K$, this morphism induces a morphism of MV-algebras $\xi_1 \colon (M \lor N)/I \rightarrow \initial$, defined by $\xi_1([x])=1$ if and only if $x \in \lnot K$. Similarly, we define a morphism of MV-algebras $\xi_2 \colon (M \lor N)/J \rightarrow \initial$ by $\xi_2([x])=1$ if and only if $x \in \lnot K$.
Let us then show that the diagram
% https://q.uiver.app/?q=WzAsNCxbMCwwLCJNIFxcbG9yIE4iXSxbMSwwLCIoTSBcXGxvciBOKS9KIl0sWzEsMSwiXFxpbml0aWFsIl0sWzAsMSwiKE0gXFxsb3IgTikvSSJdLFswLDEsInFfSiIsMCx7InN0eWxlIjp7ImhlYWQiOnsibmFtZSI6ImVwaSJ9fX1dLFsxLDIsIlxceGlfMiIsMCx7InN0eWxlIjp7ImhlYWQiOnsibmFtZSI6ImVwaSJ9fX1dLFszLDIsIlxceGlfMSIsMix7InN0eWxlIjp7ImhlYWQiOnsibmFtZSI6ImVwaSJ9fX1dLFswLDMsInFfSSIsMix7InN0eWxlIjp7ImhlYWQiOnsibmFtZSI6ImVwaSJ9fX1dXQ==
\[\begin{tikzcd}\label{diagramma IvJ not A}\tag{$\spadesuit$}
	{M \lor N} & {(M \lor N)/J} \\
	{(M \lor N)/I} & \initial
	\arrow["{q_J}", two heads, from=1-1, to=1-2]
	\arrow["{\xi_2}", two heads, from=1-2, to=2-2]
	\arrow["{\xi_1}"', two heads, from=2-1, to=2-2]
	\arrow["{q_I}"', two heads, from=1-1, to=2-1]
\end{tikzcd}\]
is a regular pushout. Consider an element $[x] \in (M \lor N)/I$ such that $[x] \in \kerr(\xi_1)$, which means that $x \in K$. Therefore, since $K = I \lor J$, there exist $i \in I$ and $j \in J$ such that $x=i \oplus j$. Now observe that $[x]=[i \oplus j]=[j]=q_I(j)$. Hence, the restriction $q_I \colon J \rightarrow \kerr(\xi_1)$ is surjective.\\
To conclude, we observe that the vertical arrows of both squares are central extensions for the structure $\Gamma_1$ if and only if $\Rad(A) \cap I \cap J = \{ 0 \}$ (where we use the fact that $\Rad(M \lor N) = (M \lor N) \cap \Rad(A)$). Applying Proposition \ref{intersezioni banali}, we can see that this condition is equivalent to requiring that $$P(A) \cap M \cap N \in \zeros.$$
Therefore, we can define the commutator between $M$ and $N$ with respect to the adjunction $S_1 \dashv i_1$ as $$[M,N]_{\CExtmv} \coloneqq P(A) \cap M \cap N;$$
$[M,N]_{\CExtmv}$ has the following property: it belongs to $\zeros$ if and only if the vertical arrows of \eqref{diagramma IvJ=A} or \eqref{diagramma IvJ not A} (depending on the join $I \lor J$) are central extensions for $\Gamma_1$. Finally, if at least one of the ideals is not proper (for example, if $I=A$), then the study of centrality reduces to analyzing the behavior of the regular epimorphism $q_J \colon A \twoheadrightarrow A/J$ with respect to the Galois structure $\Gamma$. Notably, observing that $N=K[q_J]$ and $M=A$, we have that $q_J$ is central if and only if $P(A) \cap A \cap N=P(A) \cap N \in \zeros$. Therefore, we can define $[A,N]_{\CExtmv} \coloneqq [A,N]_{s\mv}$.\\
One possible goal for future work is to further investigate the properties of this commutator.

\section*{Acknowledgement}

The author would like to thank Andrea Montoli for his valuable insights and guidance, which were fundamental in the development of this work.

\vskip 1cm
%%% BIBLIO
\bibliography{Biblio}

\begin{thebibliography}{10}

\bibitem{barr}
M.~Barr.
\newblock Exact categories.
\newblock {\em Exact categories and categories of sheaves}, pages 1--120, 1971.

\bibitem{semi}
F.~Borceux and D.~Bourn.
\newblock {\em Mal'cev, protomodular, homological and semi-abelian categories},
  volume 566.
\newblock Springer Science \& Business Media, 2004.

\bibitem{bourn1991normalization}
D.~Bourn.
\newblock Normalization equivalence, kernel equivalence and affine categories.
\newblock In {\em Category theory}, pages 43--62. Springer, 1991.

\bibitem{unitalcat}
D.~Bourn.
\newblock Mal'cev categories and fibration of pointed objects.
\newblock {\em Applied categorical structures}, 4(2):307--327, 1996.

\bibitem{bourn2003denormalized}
D.~Bourn.
\newblock The denormalized 3$\times$ 3 lemma.
\newblock {\em Journal of Pure and Applied Algebra}, 177(2):113--129, 2003.

\bibitem{varprot}
D.~Bourn and G.~Janelidze.
\newblock Characterization of protomodular varieties of universal algebras.
\newblock {\em Theory and Applications of categories}, 11(6):143--447, 2003.

\bibitem{CappellettiLgruppi}
A.~Cappelletti.
\newblock Categorical-algebraic properties of lattice-ordered groups.
\newblock {\em arXiv preprint arXiv:2310.09264}, 2023.

\bibitem{carboni1991diagram}
A.~Carboni, J.~Lambek, and M.~C. Pedicchio.
\newblock Diagram chasing in {M}al'cev categories.
\newblock {\em Journal of Pure and Applied Algebra}, 69(3):271--284, 1991.

\bibitem{chang1958algebraic}
C.~C. Chang.
\newblock Algebraic analysis of many valued logics.
\newblock {\em Transactions of the American Mathematical society},
  88(2):467--490, 1958.

\bibitem{cignoli2013algebraic}
R.~L. Cignoli, I.~M. d'Ottaviano, and D.~Mundici.
\newblock {\em Algebraic foundations of many-valued reasoning}, volume~7.
\newblock Springer Science \& Business Media, 2013.

\bibitem{everaert2010homology}
T.~Everaert and M.~Gran.
\newblock Homology of n-fold groupoids.
\newblock {\em Theory Appl. Categ}, 23(2):22--41, 2010.

\bibitem{everaert2015protoadditive}
T.~Everaert and M.~Gran.
\newblock Protoadditive functors, derived torsion theories and homology.
\newblock {\em Journal of Pure and Applied Algebra}, 219(8):3629--3676, 2015.

\bibitem{everaert2011galois}
T.~Everaert and T.~Van~der Linden.
\newblock Galois theory and commutators.
\newblock {\em Algebra universalis}, 65(2):161--177, 2011.

\bibitem{facchini2020pretorsion}
A.~Facchini and C.~A. Finocchiaro.
\newblock Pretorsion theories, stable category and preordered sets.
\newblock {\em Annali di Matematica Pura ed Applicata (1923-)},
  199(3):1073--1089, 2020.

\bibitem{facchini2021pretorsion}
A.~Facchini, C.~A. Finocchiaro, and M.~Gran.
\newblock Pretorsion theories in general categories.
\newblock {\em Journal of Pure and Applied Algebra}, 225(2):106503, 2021.

\bibitem{gran2004galois}
M.~Gran and V.~Rossi.
\newblock Galois theory and double central extensions.
\newblock {\em Homology, Homotopy and Applications}, 6(1):283--298, 2004.

\bibitem{huq}
S.~A. Huq.
\newblock Commutator, nilpotency and solvability in categories.
\newblock {\em Q. J. Math.}, 19(2):363--389, 1968.

\bibitem{idziak1984lattice}
P.~M. Idziak.
\newblock Lattice operations in bck-algebras.
\newblock {\em Mathematica Japonica}, 29:839--846, 1984.

\bibitem{janelidze1990pure}
G.~Janelidze.
\newblock Pure {G}alois theory in categories.
\newblock {\em Journal of Algebra}, 132(2):270--286, 1990.

\bibitem{semiabel}
G.~Janelidze, L.~M{\'a}rki, and W.~Tholen.
\newblock Semi-abelian categories.
\newblock {\em Journal of Pure and Applied Algebra}, 168(2-3):367--386, 2002.

\bibitem{librogruppiret}
V.M. Kopytov and N.Y. Medvedev.
\newblock {\em The theory of lattice-ordered groups}, volume 307.
\newblock Springer Science \& Business Media, 2013.

\bibitem{lapenta2022relative}
S.~Lapenta, G.~Metere, and L.~Spada.
\newblock Relative ideals in homological categories, with an application to
  mv-algebras.
\newblock {\em arXiv preprint arXiv:2208.12597}, 2022.

\bibitem{mundici1986interpretation}
D.~Mundici.
\newblock Interpretation of {AF} {C}\text{*}-algebras in {\l}ukasiewicz
  sentential calculus.
\newblock {\em Journal of Functional Analysis}, 65(1):15--63, 1986.

\bibitem{aritmetica}
M.~C. Pedicchio.
\newblock Arithmetical categories and commutator theory.
\newblock {\em Applied Categorical Structures}, 4:297--305, 1996.

\bibitem{pixley1963distributivity}
A.~F. Pixley.
\newblock Distributivity and permutability of congruence relations in
  equational classes of algebras.
\newblock {\em Proceedings of the American Mathematical Society},
  14(1):105--109, 1963.

\bibitem{salas1980estudio}
A.~J.~R. Salas.
\newblock {\em Un estudio algebraico de los c{\'a}lculos proposicionales de
  Lukasiewicz}.
\newblock PhD thesis, Universitat de Barcelona, 1980.

\end{thebibliography}
\bibliographystyle{plain}
\end{document}